\documentclass[12pt,a4paper]{article}
 \usepackage[latin1]{inputenc}
 \usepackage{verbatim}
 \usepackage[dvips]{color}      
 \usepackage{flafter}
\newcounter{myenumii}
\newcounter{anhz}

\refstepcounter{anhz}

\renewcommand{\theequation}{\arabic{section}.\arabic{equation}}

\newtheorem{lemma}{Lemma}[section]
\newtheorem{theorem}[lemma]{Theorem}
\newtheorem{corollary}[lemma]{Corollary}
\newtheorem{prop}[lemma]{Proposition}

\newtheorem{remark}[lemma]{Remark}

\newcommand{\N}{{\ensuremath{\mathrm{I}\!\mathrm{N}}}}
\newcommand{\R}{{\ensuremath{\mathrm{I}\!\mathrm{R}}}}

\definecolor{orange}{rgb}{1,0.5,0} 
\definecolor{gelb}{rgb}{1,1,0}
\definecolor{blaugrau}{rgb}{0,0.5,0.5}
\definecolor{dgreen}{rgb}{0,.8,0.4}
\definecolor{hellblau}{rgb}{0,0.8,1}
\definecolor{hellblau2}{rgb}{0,0.5,1}
\definecolor{violet}{rgb}{0.5,0,0.5}
\definecolor{lila}{rgb}{1,0,1}
\definecolor{grau}{rgb}{0.3,0.3,0.3}
\definecolor{hellgrau}{rgb}{0.85,0.85,0.85}
\definecolor{dlila}{rgb}{0.8,0.0,0.4}
\definecolor{braun}{rgb}{0.9,0.6,0}
\definecolor{braunlila}{rgb}{0.4,0.1,0.2}

\newcommand{\diver}{\mbox{\rm div}}
\newcommand{\supp}{\mbox{\rm supp}}

\newcommand{\qed}{\hfill\framebox[0.2cm]}
\newcommand{\proof}{{\sc Proof}}

\newcommand{\INT}{\!\int\!\!\!\!\int\!\!}

\newcommand{\fer}[1]{(\ref{#1})}

\begin{document}

\title{ The Wigner-Poisson-Fokker-Planck system: global-in-time solution and dispersive
effects  }
\author{Anton Arnold, \and
        Elidon Dhamo,\and
        Chiara Manzini}
\date{}
\maketitle

\medskip


\abstract{This paper is concerned with the
Wigner-Poisson-Fokker-Planck system, a kinetic evolution equation
for an open quantum system with a non-linear Hartree potential.
Existence, uniqueness and regularity of global solutions to the
Cauchy problem in 3 dimensions are established. The analysis is
carried out in a weighted $L^2$--space, such that the linear
quantum Fokker-Planck operator generates a dissipative semigroup.
The non-linear potential can be controled by using the parabolic
regularization of the system.

The main technical difficulty for establishing global-in-time
solutions is to derive a-priori estimates on the electric field:
Inspired by a strategy for the classical Vlasov-Fokker-Planck
equation, we exploit dispersive effects of the free
transport operator. As a ``by-product'' we also derive a new a-priori
estimate on the field in the Wigner-Poisson equation.}\\

\newpage
\setlength{\parindent}{0cm}

\section{Introduction}
\setcounter{equation}{0}

The goal of this paper is to prove the existence and uniqueness of global-in-time solutions to
the coupled Wigner-Poisson-Fokker-Planck (WPFP) system  in three dimensions.  This kinetic
equation is the quantum mechanical analogue of the classical Vlasov-Poisson-Fokker-Planck (VPFP)
system, which models the diffusive transport of charged particles (in plasmas e.g.).

Wigner functions provide a kinetic description of quantum
mechanics (cf. \cite{Wig32}) and have recently become a valuable
modeling and simulation tool in fields like semiconductor device
modeling (cf. \cite{MaRiSc90} and references therein), quantum
Brownian motion, and quantum optics (\cite{CaLe83,Dio93a}). The
real-valued Wigner function $w(x,v,t)$ is a probabilistic
quasi-distribution function in the position-velocity $(x,v)$ phase
space for the considered quantum system at time $t$.

Its temporal evolution is governed by the Wigner-Fokker-Planck
(WFP) equation:
\begin{equation}
\label{eq:1.1} w_{t} + v\cdot\nabla_xw - \Theta[V]w \;\: = \;\:
\beta\diver_v(vw)+\sigma \Delta_vw+2\gamma \diver_v(\nabla_xw) + \alpha \Delta_xw,\quad
t>0,
\end{equation}
on the phase space $x\in\R^3,\; v\in\R^3,$ with the initial condition
$$
w(x,v,t=0) \;\: = \;\: w_0(x,v).
$$
With a vanishing right hand side, equation \fer{eq:1.1} would be
the (diffusion-free) Wigner equation. It describes the reversible
evolution of a quantum system under the action of a (possibly
time-dependent) electrostatic potential $V=V(x,t)$. The potential
effect enters in the equation via the pseudo-differential operator
$\Theta[V]$:
\begin{eqnarray}
\label{eq:1.4} (\Theta[V]w)(x,v,t) & = &
i[V(t,x+\frac{1}{2i}\nabla_{v}) -
V(t,x-\frac{1}{2i}\nabla_{v})]w(x,v,t)  \nonumber\\ & = &
\frac{i}{(2\pi)^{3/2}}\int_{\R^3}\delta V(x,\eta,t){\cal
F}_{v\rightarrow\eta} w(x,\eta,t)e^{iv\cdot\eta}\,d\eta
\nonumber\\ & = & \frac{i}{(2\pi)^3} \int_{\R^3}\int_{\R^3}\delta
V(x,\eta,t)w(x,v^{,},t)e^{i(v-v^{,})\cdot\eta}\,dv^{,}\,d\eta,
\end{eqnarray}
where $\delta V(x,\eta,t) = V(x+\frac{\eta}{2},t) -
V(x-\frac{\eta}{2},t)$ and ${\cal F}_{v\rightarrow\eta} w$ denotes
the Fourier transform of $w$ with respect to $v$:
\begin{eqnarray*}
{\cal F}_{v\rightarrow\eta} w(t,x,\eta) \;\:=\;\: \frac{1}{(2\pi)^{3/2}}
\int_{\R^3}w(t,x,v^{,})e^{-iv^{,}\cdot\eta}\,dv^{,}.
\end{eqnarray*}
For simplicity of the notation we have set here the Planck
constant, particle mass and charge equal to unity.\\

The right hand side of \fer{eq:1.1} is a Fokker-Planck type model
for the non-reversible interaction of this quantum system with an
environment, e.g.\ the interaction of electrons with a phonon bath
(cf.\ \cite{CaLe83,CEFM00} for derivations from reversible quantum
systems, and \cite{GGKS93,Str86} for applications in quantum
transport). In \fer{eq:1.1}, $\beta\geq0$ is the friction
parameter and the parameters $\alpha,\;\gamma \geq 0$, $\sigma >
0$ constitute the phase-space diffusion matrix of the system. In
the kinetic Fokker-Planck equation of classical mechanics (cf.
\cite{Ris84,CSV96}) one would have $\alpha=\gamma=0$. For the WFP
equation \fer{eq:1.1} we have to assume
$$ \left(\begin{array}{cc} \alpha & \gamma +
\frac{i}{4}\beta \\ \gamma - \frac{i}{4}\beta & \sigma
\end{array}\right) \;\:\ge \;\: 0,
$$ which guarantees that the system is {\it quantum mechanically
correct}. More precisely, it guarantees that the corresponding von
Neumann equation is in Lindblad form and that the density matrix
of the quantum system stays a positive operator under temporal
evolution (see \cite{ALMS00} for details).
In the sequel we shall therefore assume
\begin{equation}
\label{coeff} \alpha\sigma \;\:\ge\;\: \gamma^2 + \frac{\beta^2}{16}\quad \mbox{and}\quad
\alpha\sigma \;\:>\;\: \gamma^2.
\end{equation}
Hence, the principle part of the Fokker-Planck term is uniformly elliptic. This makes the
present work complementary to \cite{ALMS00}, where the friction-free, hypoelliptic case (with
$\alpha=\beta=\gamma=0$) was analyzed.\\

The WFP equation \fer{eq:1.1} is self-consistently coupled with
the Poisson equation for the (real-valued) potential
$V=V[w](x,t)$:
\begin{equation}
-\Delta V  \;\:=\;\:  n[w],\quad x\in\R^3,\quad t>0,
\end{equation}
with the \emph{particle density}
\begin{equation}
\label{eq:1.7} n[w](x,t) \;\::=\;\: \int_{\R^3} w(x,v,t)\,dv.
\end{equation}
This potential models the repulsive Coulomb interaction within the considered
particle system in a mean-field description. \\

The main analytical challenge for tackling Wigner-Poisson systems is the proper definition of
$n[w]$ in appropriate $L^p$ spaces. Due to the definition of the operator $\Theta$ in Fourier
space, $w\in L^2(\R^3_x\times \R^3_v)$ is the natural set-up. Without further assumptions, of
course, this does not justify to define $n[w]$. We shall now summarize the existing literature
of this field and the typical strategies to overcome the above problem:

a) The standard approach for the \emph{Wigner-Poisson} equation is to reformulate it as a
Schr\"o\-din\-ger-Poisson system, where the particle density then appears in $L^1$ (cf.\
\cite{BrMa,Ca1} for the 3D--whole space case).

b) In \emph{one spatial dimension with periodic boundary conditions} in $x$ the
Wigner-Poisson system (and WPFP) can be dealt with directly on the kinetic level. For $w$ in
a weighted $L^2$-space, the nonlinear term $\Theta[V]w$ is then bounded and locally
Lipschitz \cite{ArRi96,ArCaDh}. The same strategy was also used in \cite{Ma03} for the
Wigner-Poisson system on a bounded (spatial) domain in three dimensions (local-in-time
solution).

c) By adapting $L^1$-techniques from the classical
Vlasov-Fokker-Planck equation, the \emph{3D
Wigner-Poisson-Fokker-Planck} system was analyzed in \cite{ALMS00}
(local-in-time solution for the friction-free problem) and
\cite{CaLoNi04} (global-in-time solution). The latter paper,
however, is not a purely kinetic analysis as it requires to assume
the positivity of the underlying density matrix. In both cases the
dissipative structure of the system allows to control $n[w]$.

d) In \cite{Arn95a,ArSp} the 3D Wigner-Poisson and WPFP systems
were reformulated as von-Neumann equations for the quantum
mechanical \emph{density matrix}. This implies $n\in L^1(\R^3)$.
While this approach is the most natural, both physically and in
its mathematical structure, it is restricted to whole space cases.
Extensions to initial-boundary value problems (as needed for
practical applications and numerical analysis) seem
unfeasible.

e) For the classical \emph{Vlasov-Poisson-Fokker-Planck} equation there exists
a vast body of mathematical literatur from the 1990's (cf.\
\cite{Bou93,Bou95,CSV96,Ca,Ca98}), and many of those tools will be closely related to
the present work.\\

In spite of the various existing well-posedness results for the WPFP problem,
there is a need for a purely kinetic analysis, and this is our goal here. Such an approach could possibly allow for an extension to boundary-value problems in the Wigner framework (where the positivity of the related density matrix  is a touchy question). 

Mathematically we shall develop the following new tools and
estimates that could be important also for other quantum kinetic
applications: In all of the existing literature on Wigner-Poisson
problems (except \cite{St91}) the potential $V$ is bounded, which
makes it easy to estimate the operator $\Theta[V]$ in $L^2$. Our framework for the local 
in time analysis does not yield a bounded potential. However, the
operator $\Theta$ only involves $\delta V$, a potential
difference, which has better decay properties at infinity. This
observation gives rise to new estimates that are crucial for our
local-in-time analysis.

In order to establish global-in-time solutions we shall extend dispersive tools
of Lions, Perthame and Castella (cf.\ \cite{LiPe91,Pe96, CaPe96} for applications to
classical kinetic equation) to the WP and WPFP systems.  The fact that the
Wigner function $w$ also takes negative values gives rise to an important
difference between classical and quantum kinetic problems:  In the latter case,
the conservation of mass and energy or pseudo-conformal laws do \emph{not}
provide useful a-priori estimates on $w$. We shall hence assume that the
initial state lies in a weighted $L^2$-space, but we shall \emph{not} require
that our system has finite mass or finite kinetic energy. Since the energy
balance will not be used, this also implies that the sign of the interaction
potential does \emph{not} play a role in our analysis.\\

This paper is organized as follows: In Section 2 we introduce a
weighted $L^2$-space for the Wigner function $w$ that allows to
define $n[w]$ and the nonlinear term $\Theta[V]w$. In \S3 we
obtain a local-in-time, mild solution for WPFP using a fixed point
argument and the parabolic regularization of the Fokker-Planck
term. In \S4 we establish a-priori estimates to obtain
global-in-time solutions. The key point is to derive first
$L^p$-bounds for the electric field $\nabla V$ by exploiting
dispersive effects of the free kinetic transport. ``Bootstraping''
then yields estimates on the Wigner function in a weighted
$L^2$-space. Finally, we give regularity results on the solution.
The technical proofs of several lemmata are defered to the
Appendix.

\bigskip

\noindent {\bf Acknowledgement}\\
This work was supported by the IHP--Network HPRN-CT-2002-00282
from the European Union. The first author acknowledges support
from the DFG under Grant-No. AR 277/3-2. The third author was partially
supported by the GNFM Research Project ``Mathematical Models for
Microelectronics'', 2004 and MIUR-COFIN ``Mathematical
problems of Kinetic Theories'', 2002.


\section{The functional setting}
\setcounter{equation}{0}
In this section we shall discuss the functional analytic preliminaries for
studying the non-linear problem \fer{eq:1.1}-\fer{eq:1.7}. First we shall
introduce an appropriate ``state space'' for the Wigner function $w$ which
allows to ``control'' the particle density $n[w]$ and the selfconsistent
potential $V[w]$. Next, we shall discuss the linear Wigner-Fokker-Planck
equation and the dissipativity of its (evolution) generator $A$.

\subsection{State space and selfconsistent potential}
\label{subsection:state_space}

Let us introduce the following weighted (real valued) $L^2$-space
\begin{equation}
X \;\::=\;\: L^{2}(\R^6; (1 + |v|^2)^2\,dx\,dv),
\end{equation}
endowed with the scalar product
\begin{equation}
<u,w>_{X} \;\:=\;\: \int_{\R^3}\!\int_{\R^3}\!u\,w(1+|v|^2)^2\,dx\,dv.
\end{equation}
The following proposition motivates the choice of $X$ as the
state space for our analysis.
\begin{prop}
\label{prop:density} For all $w\in X,$ the function $n[w]$ defined
by $n[w](x):=\int\!w(x,v)\,dv,$\\$x\in \R^3,$ belongs to
$L^{2}(\R^3)$ and satisfies
\begin{equation}
\label{eq:density1}
\Vert n[w]\Vert_{L^{2}(\R^3)} \;\:\leq\;\: C\Vert w\Vert_{X},
\end{equation}
with a constant $C$ independent of $w.$
\end{prop}
Here and in the sequel $C$ shall denote generic, but not necessarily equal, constants.\\
\proof.
By using H\"older inequality in the $v$-integral, we get
\begin{eqnarray*}
\Vert n[w]\Vert_{L^{2}(\R^3)}^2 
 & \leq &   \!\!\int\!\!\left(\int\!\!|w(x,v)|^2(1 + |v|^2)^2 \,dv\right) \left(\int\!\!\frac{dv}{(1 + |v|^2)^2} \right)dx 
\;\: = \;\: C\Vert w\Vert_{X}^{2}.
\end{eqnarray*}
$\quad $ \hfill\framebox[0.2cm]\\

\begin{remark}
The choice of the $|v|^2$ weight was already seen to be
convenient to control the $L^2$-norm of the density on a
bounded domain of $\R_x^3$ (cf.\ \cite{Ma03}). \\
The subsequent analysis would hold also by including a symmetric weight  in the $x$-variable (i.e. for $w\in L^{2}(\R^6; (1 + |x|^2+|v|^2)^2\,dx\,dv)$ ), which would yield a $L^p$-bound with $p\in(3/2,2]$ for the density. 
On the other hand, Lemmata \ref{lemma:dissipA}, \ref{lemma:density} would prevent us from introducing a non-symmetric weight in $x$.
\end{remark}

In this framework
the following estimates for the self-consistent potential hold.
\begin{prop}\label{prop:potentialI} For all $w\in X,$ the (Newton potential) solution $V=V[w]$ of the equation
$-\Delta_x V[w]=n[w],\,x\in \R^3,$ satisfies
\begin{equation}
\label{eq:potentialI}
\Vert\nabla V[w]\Vert_{L^{6}(\R^3)} \;\:\leq\;\: C \Vert
n[w]\Vert_{L^{2}(\R^3)}.
\end{equation}
\end{prop}
\proof.
Since $V=-\frac{1}{4\pi|x|}*n,$ we have $\nabla
V=\frac{x}{4\pi|x|^3}*n,$ and the estimate follows from the
generalized Young inequality.\hfill\framebox[0.2cm]\\
\begin{remark}
Note that $n\in L^2(\R^3)$ does not yield (via the generalized Young inequality) a control
of $V$ in any $L^r-$space  (even  $n\in L^p(\R^3)$ with $p\in(3/2,2]$ would not ``help''). 
However, the operator $\Theta[V]$ involves only the function $\delta
V,$ which is slightly ``better behaved''.\\
We anticipate that we will later recover some information on the potential $V$ via new \emph{a-priori} estimates on the electric field $\nabla V[w]$ (see Corollary \ref{cor:WPFP-potential}).
\end{remark}
Omitting the time-dependence we have
\begin{eqnarray*}
\delta V(x,\eta) &=& V(x+\frac{\eta}{2}) - V(x-\frac{\eta}{2}) \;\:=\;\:
\frac{1}{4\pi} \int\limits_{\R^3} \frac{n[w](x-\frac{\eta}{2}-\xi)-n[w](x+\frac{\eta}{2}-\xi)}{|\xi|}
\,d\xi \\ &=& \frac{1}{4\pi}\int\limits_{\R^3} f(y;\eta) \,n[w](x-y)\,dy,
\end{eqnarray*}
with the ``dipole-kernel'' $f(y;\eta):=\left(
\frac{1}{|y-\frac{\eta}{2}|}
-\frac{1}{|y+\frac{\eta}{2}|}\right).$
\begin{prop}
\label{prop:potentialII} For all $w\in X$ and fixed $\eta \in \R^3,$ we have $\delta
V[w](.,\eta) \in L^q(\R_x^3), \; 6< q\leq \infty. $ Moreover
\begin{equation}
\label{eq:potentialII}
\Vert \delta V[w](.,\eta)\Vert_{L^\infty(\R_x^3)} \;\:\leq \;\: C|\eta|^{1/2}\Vert n[w]\Vert_{L^2(\R^3)}.
\end{equation}
\end{prop}
\proof.
By using the triangle inequality, $$
    |f(y;\eta)| \;\:\leq\;\:
    \frac{|\eta|}{|y-\frac{\eta}{2}||y+\frac{\eta}{2}|}\,,
$$ and the transformation $y=|\eta|x,$ we estimate for $3/2<p< 3$
$$ \Vert f(\ldotp ;\eta)\Vert_{L^p(\R^3)}^{p} \;\:=\;\:
|\eta|^{3-p}\int\limits_{\R^3}
\frac{dx}{\left(|x-\frac{e}{2}||x+\frac{e}{2}|\right)^p} \;\:<\;\:\infty,
$$ where $e\in \R^3$ is some unit vector (due to the rotational symmetry of $\Vert f(.;\eta)\Vert_{L^p(\R^3)}^{p}$ with respect to $\eta$).
Young inequality then gives $\delta V(.,\eta) \in L^q(\R^3),$
$6<q\leq\infty,$ and the assertion holds.
\qed\\

\bigskip
In most of the literature the {\it Wigner operator} $\Theta$ is
defined on $L^2(\R_{v}^{d})$ for bounded potentials $V,$ cf.
\cite{MaRi89,MaBa03,ArCaDh}. For our nonlinear problem
\fer{eq:1.1}-\fer{eq:1.7}, however, $V\in L^\infty(\R^3)$ does
{\it not} hold. As a compensation we shall hence exploit the
additional regularity of the Wigner function to define the
quadratic term $\Theta[V[w]]w$ (cf. Prop. 2.8 in \cite{Ma03} for a
similar strategy).
\begin{prop}
\label{prop:pert} Let $u\in X$ and $\nabla_v u \in X$ be given.
Then, the linear operator $$ z\longmapsto \Theta[V[z]]u,$$ with
the function $V[z]=-\frac{1}{4\pi|x|}*n[z]$, is bounded from the
space $X$ into itself and satisfies
\begin{equation}
\label{eq:pert1} \Vert \Theta[V[z]]u\Vert_{X}\;\:\leq\;\:
C\{\,\Vert u\Vert_{X}+\Vert \nabla_v u\Vert_{X}\}\,\Vert
z\Vert_{X}, \quad \forall z\in X.
\end{equation}
\end{prop}
\proof. To estimate $\Vert \Theta[V[z]]u\Vert_X$ we shall consider
separately the two terms of the equivalent norm

\begin{equation}
\label{eq:X-tilde} \Vert u \Vert_{\tilde{X}}^{2} \;\::=\;\: \Vert
u\Vert_{2}^{2} + \sum_{i=1}^{3} \Vert v_i^2u \Vert_2^2.
\end{equation}

First, by denoting $\hat{u}:={\cal F}_{v\rightarrow\eta} u,$ we
get
\begin{eqnarray}
\Vert\Theta[V[z]]u\Vert_2^2 &=& \!\int\!\!\!\int \!\!|\delta(V[z])(x,\eta)\hat{u}(x,\eta)|^2dx\,d\eta
\;\:\leq\:\;\!\int\!\!\!\int \!\!\Vert \delta(V[z])(\ldotp,\eta)\Vert_{\infty}^2|\hat{u}(x,\eta)|^2 d\eta\,dx\nonumber\\
\label{eq:pertL2}
&\leq&C\Vert
z\Vert_{X}^2\!\int\!\!\!\int\!\!\left(|\eta|^{1/2}|\hat{u}(x,\eta)|\right)^2 d\eta\,dx
\;\:\leq\:\; C\Vert z\Vert_{X}^2 \left(\Vert u\Vert_{2}^2+\Vert\nabla_v u\Vert_{2}^2\right),
\end{eqnarray}
by applying first the Plancherel Theorem, then H\"older's inequality in
the $x$ variable, the estimates (\ref{eq:density1}), (\ref{eq:potentialII}) for the function $\delta V[z],$ and
finally, Young inequality and the Plancherel Theorem to the last
integral.\\
For the second term of $\Vert\Theta[V[z]]u\Vert_{\tilde X}$ we
shall use
\begin{equation}
\label{eq:2.9}
v_i^2\Theta[V]w(x,v) \;\:=\;\: \frac{1}{4}\Theta[\partial_i^2V]w(x,v) +
\Omega[\partial_i V](v_i w)(x,v)+\Theta[V]v_i^2w(x,v),
\end{equation}
with the pseudo-differential operator
\begin{equation}
\label{eq:defOmega}
\Omega[V]\;\::=\;\:i(\delta_+V)\left(x,\frac{\nabla_v}{i}\right),\qquad
(\delta_+ V)(x,\eta)\;\::=\;\:V\! \left( x +
\frac{\eta}{2} \right) +V\!\left(x - \frac{\eta}{2}\right).
\end{equation}
Here and in the sequel we use the abreviation $\partial_i:=\partial_{x_i}$.
\fer{eq:2.9} is now estimated:
\begin{equation}
\label{eq:pertVweight}
\Vert v_i^2\Theta[V[z]]u\Vert_2 \;\:\leq\;\: \frac{1}{4}\Vert \delta(\partial_i^2
V[z]) \hat{u}\Vert_2+\Vert\delta_+(\partial_i
V[z])\partial_{\eta_i}\hat{u}\Vert_2+ \Vert\delta
V[z]\partial_{\eta_i}^2\hat{u}\Vert_2
\end{equation}
The first two terms of (\ref{eq:pertVweight}) can be estimated as follows:
\begin{eqnarray*}
\Vert \delta(\partial_i^2 V[z]) \hat{u}\Vert_{L^2(\R^6)}&\leq& 2\Vert
\partial_i^2 V[z]\Vert_{L^2(\R_x^3)}\Vert \hat{u}\Vert_{L^2(\R_x^3;L^{\infty}(\R_\eta^3))}\\
&\leq& C\Vert z
\Vert_X\Vert (1+|v|^2)u\Vert_{L^2(\R^6)},
\end{eqnarray*}
by applying H\"older's inequality, (\ref{eq:density1}) and the Sobolev imbedding 
$\hat{u}(x,\ldotp)\in H^2(\R_{\eta}^3)\hookrightarrow
L^{\infty}(\R_{\eta}^3).$
\begin{eqnarray}
\label{eq:deltaplus}
\Vert\delta_+(\partial_i
V[z])\partial_{\eta_i}\hat{u}\Vert_{L^2(\R^6)}&\leq& C \Vert \partial_{i}
V[z]\Vert_{L^4(\R_x^3)}\Vert\partial_{\eta_i}\hat{u}\Vert_{L^2(\R_x^3;L^4(\R_{\eta}^3))}
\nonumber\\
&\leq& C\Vert z\Vert_{X} \Vert  (1+v_i^2)u\Vert_2,
\end{eqnarray}
by the Sobolev imbedding and $\nabla_{\eta}\hat{u}(x,\ldotp)\in H^1(\R_{\eta}^3)\hookrightarrow L^4(\R_{\eta}^3),$
and by estimate (\ref{eq:potentialI}) for $\nabla V[z]$ and (\ref{eq:density1}). For the last
term of (\ref{eq:pertVweight}) we estimate as in (\ref{eq:pertL2}):
\begin{eqnarray*}
\Vert\delta V[z]\partial_{\eta_i}^2\hat{u}\Vert_2^2
&\leq&\int\!\!\!\int \!\!\Vert \delta V[z](\ldotp,\eta)\Vert_{\infty}^2|\partial_{\eta_i}^2\hat{u}(x,\eta)|^2d\eta\,dx\\
&\leq& C\Vert z\Vert_{X}^2\int\!\!\!\int\!\!\left(|\eta|^{1/2}\partial_{\eta_i}^2\hat{u}(x,\eta)\right)^2d\eta\,dx\\
&\leq& C\Vert z\Vert_{X}^2\left(\Vert
\partial_{\eta_i}^2\hat{u}\Vert_2^2 +\Vert\eta
\partial_{\eta_i}^2\hat{u}\Vert_2^2\right)\\ &\leq& C\Vert
z\Vert_{X}^2\left(\Vert\partial_{\eta_i}^2\hat{u}\Vert_2^2 +
\Vert\partial_{\eta_i}^2(\eta \hat{u})\Vert_2^2 + \Vert
\partial_{\eta_i}\hat{u}\Vert_2^2\right)\\
&\leq& C\Vert z\Vert_{X}^2\left(\Vert (1+v_i^2)u\Vert_2^2+ \Vert
v_i^2\nabla_v u\Vert_2^2 \right),
\end{eqnarray*}
by interpolation and integration by parts. \\
This concludes the proof of estimate (\ref{eq:pert1}). \hfill\framebox[0.2cm]\\

\begin{remark}
The previous proposition shows that the bilinear map
$$(z,u)\longmapsto\Theta[V[z]]u$$
is well-defined for all $z,u\in X,$ subject to $\nabla_v u \in X.$
The unusual feature of the above proposition is the boundedness of this map
with respect to the function $z$ appearing in the self-consistent
potential $V[z].$ This is in contrast to most of the existing
literature (\cite{ArCaDh,MaBa03,MaRi89}), where the boundedness of
the pseudo-differential operator $\Theta[V[z]]$ (with $z$ fixed)
is used. However, this can only hold for bounded potentials $V$.
\end{remark}

\subsection{Dissipativity of the linear equation}
\label{subsection:dissipativity}
In our subsequent analysis we shall first consider the linear
Wigner-Fokker-Planck equation, i.e. equation \fer{eq:1.1} with $V\equiv0$. The
generator of this evolution problem is the unbounded linear operator $A : D(A)
\longrightarrow X,$
\begin{equation}
\label{def:A}
Au \;\: := \;\: -v\cdot\nabla_{x}u + \beta\diver_{v}(vu) + \sigma\Delta_{v}u
+ 2\gamma\diver_v(\nabla_x u)+ \alpha\Delta_{x}u,
\end{equation}
defined on
\begin{eqnarray}
D(A)\;\: = \;\: \lbrace u\in X \;|\; v\cdot\nabla_{x}u, v\cdot\nabla_vu,\Delta_{v}u,\diver_v\nabla_x u, \Delta_{x}u\in X \rbrace.
\end{eqnarray}
Clearly, $C^{\infty}_{0}(\R^6)\subset D(A).$ Hence, $D(A)$ is
dense in $X.$ Next we study whether the operator $A$ is dissipative
on the (real) Hilbert space $\tilde X$, i.e.\ if
\begin{equation}
<Au,u>_{\tilde X} \;\:\leq\;\: 0, \quad\forall \, u\in D(A)
\end{equation}
holds.

\begin{lemma}
\label{lemma:dissipA}
Let the coefficients of the operator $A$ satisfy $\alpha\sigma\geq\gamma^2.$ Then $A- \kappa I$ with
\begin{equation}
\label{def:kappa}
\kappa \;\;:=\;\;  \frac{3}{2}\beta  + 9\sigma
\end{equation}
 is dissipative in $\tilde{X}$.
\end{lemma}

The proof is lengthy but straigthforward and deferred to the Appendix.

\vspace*{0.3cm}
By Theorem 1.4.5b of \cite{Pazy} its closure, $\overline{A -
\kappa I} = \overline{A} - \kappa I$ is also dissipative. \\ A
straightforward calculation using integrations by parts yields
$$
<Au,w>_{\tilde{X}} \;\:=\;\: <u,A^{*}_{1} w>_{\tilde{X}} \,+\,
<u,A^{*}_{2} w>_{\tilde{X}}, \quad \forall\, u,w\in D(A),
$$
with
\begin{eqnarray*}
A^{*}_{1}w & = & v\cdot\nabla_x w - \beta v\cdot\nabla_v w
+ \sigma \Delta_v w + 2\gamma\diver_v(\nabla_x w)+ \alpha\Delta_{x}w, \\
<u,A^{*}_{2} w>_{\tilde{X}} & = & \sum_{i=1}^{3}\left(  - \frac{4}{3}\beta \INT v^{4}_{i}wu + \frac{8}{3}\sigma\INT v^{3}_{i}w_{v_i}u + \frac{12}{3}\sigma\INT v^{2}_{i}wu + \frac{8}{3}\gamma\INT v^{3}_{i}w_{x_i}u\right).
\end{eqnarray*}

\vspace*{0.3cm}
Hence, $A^*\!\!\mid_{D(A)}$ -- the restriction of the adjoint of the
operator $A$ to $D(A)$ -- is given by $A^*
\,=\,A^{*}_{1}+A^{*}_{2}.$ $A^*$ is densly defined on $D(A^*)
\supseteq D(A),$ and hence $A$ is a closable operator (cf. Theorem
VIII.1.b of \cite{RS1}). Its closure $\overline{A}$ satisfies
$(\overline{A})^* = A^*$ (cf. \cite{RS1}, Theorem VIII.1.c). \\
Since $<A^*u,u>$ = $<Au,u>$ the following lemma on the dissipativity of the operator $A^*$ restricted to
$D(A)$ holds.
\begin{lemma}
\label{lemma:dissipA*}
Let the coefficients of the operator $A$ satisfy $\alpha\sigma\geq\gamma^2.$
Then $A^*\!\!\mid_{D(A)} - \kappa I$ is dissipative (with $\kappa$ as in (\ref{def:kappa})).
\end{lemma}

Next we consider the dissipativity of this operator on its proper domain $D(A^*),$ which, however,
is not known explicitly. To this end we shall use the following technical lemma, which, for a matter of generality, is stated in the space with the symmetric $x,v$-weight. The proof is defered to the
appendix: the arguments are inspired by \cite{ArCaDh}, \cite{ArSp}, but there are also similar results for FP-type
operators in \cite{HeNi03,HeNi04}, e.g.

\begin{lemma}
\label{lemma:density}
Let $P=p(x,v,\nabla_x,\nabla_v)$ where $p$ is a quadratic polynomial and
$$
 D(P) \;\: := \;\: C^{\infty}_{0}(\R^6) \subset Z\;\: := \;\:L^{2}(\R^6; (1 + |x|^2+|v|^2)^2\,dx\,dv).
$$
Then $\bar{P}$ is the maximum extension of $P$ in the sense that
$$
 D(\bar{P}) \;\: := \:\; \left\{ u\in Z \:| \mbox{ the distribution } Pu\in Z\right\}.
$$
\end{lemma}

We now apply Lemma \ref{lemma:density} to $P=A^* - \kappa I,$ which is dissipative on $D(P)\subset D(A).$ Since $A^*$ is
closed, we have $D(A^*) = D(\overline{P})=\{ u\in X \,|\, A^* u\in X \}$ and $A^* - \kappa I$ is dissipative
on all of $D(A^*).$

Applying Corollary 1.4.4 of \cite{Pazy} to $\overline{A} - \kappa I$ (with $(\overline{A})^* = A^*$), then implies that
$\overline{A} - \kappa I$ generates a $C_0$ semigroup of contractions on $X,$ and the $C_0$ semigroup generated by
$\overline{A}$ satisfies
$$
  \Vert e^{t\overline{A}}u \Vert_{\tilde{X}} \;\:\leq\;\: e^{\kappa t}\Vert u\Vert_{\tilde{X}}, \quad u\in X, \;\;t\geq 0.
$$
Since $\Vert .\Vert_{X}$ and $\Vert . \Vert_{\tilde{X}}$ are equivalent norms in $X$ with
$$
\Vert u\Vert_{\tilde{X}} \;\:\leq\;\: \Vert u\Vert_{X} \;\:\leq\;\: 4\Vert u\Vert_{\tilde{X}},
$$
we have
\begin{equation}
\label{eq:estimate}
\Vert e^{t\overline{A}}u \Vert_{X} \;\:\leq\;\: 4 e^{\kappa t}\Vert u\Vert_{X}, \quad u\in X, \;\;t\geq 0.
\end{equation}

\section{Existence of the local-in-time solution}
\setcounter{equation}{0}

In this section we shall use a contractive fixed point map to establish a local solution
of the WPFP system. To this end the parabolic regularization of the linear WFP equation will be
crucial to define the self-consistent potential term.

\subsection{The linear equation}
First let us consider the linear equation
\begin{equation}
\label{eq:linEq} w_{t} \;\:=\;\: \overline{A}w(t),\quad t>0,\qquad w(t=0)\;\:=\;\:w_0 \in X.
\end{equation}
By the discussion in Subsection \ref{subsection:dissipativity},
its unique solution $w(t)=e^{t\overline{A}}w_0$ satisfies
\begin{equation}
\label{eq:decayI}
\Vert w(t)\Vert_{X}\;\:\leq\;\: 4 e^{\kappa t}\Vert w_0\Vert_{X},\quad \forall\,t\geq 0.
\end{equation}

Actually, the solution of the equation can be expressed as
\begin{equation}
\label{eq:solution} w(x,v,t) \;\:=\;\: \INT
w_0(x_0,v_0)G(t,x,v,x_0,v_0)\,dx_0\,dv_0,\quad \forall\, (x,v)\in \R^6,
\end{equation}
where the Green's function $G$ satisfies (in a weak sense) the equation \fer{eq:linEq}
and the initial condition
$$
  \lim_{t\rightarrow 0}G(t,x,v, x_0,v_0)\;\:=\;\:\delta(x-x_0,v-v_0),
$$
for any fixed $(x_0,v_0)\in \R^6$
(cf.\ Def.\ 2.1 and Prop.\ 3.1 in \cite{SCDM02}). \\
The Green's function reads
\begin{eqnarray}
\label{eq:Green}G(t,x,v,x_0,v_0) & = &  e^{3\beta t} F(t,
X_{-t}(x,v)-x_0, \dot{X}_{-t}(x,v)-v_0),
\end{eqnarray}
with $$ F(t,x,v)\;\: =\:\;\frac{ 1}{(2\pi)^3\left(4\lambda(t)\nu(t)
-\mu^2(t)\right)^{3/2}}\cdot\exp\left\{ -\frac{\nu(t)|x|^2 +
\lambda(t)|v|^2 + \mu(t)(x\cdot v)}{4\lambda(t)\nu(t)
-\mu^2(t)}\right\}.$$ The characteristic flow
$\Phi_t(x,v) = [X_t(x,v),\dot{X}_{t}(x,v)]$ of the first order part of (\ref{def:A}), is
given for $\beta>0$ by
\begin{eqnarray*}
X_{t}(x,v) & = &
x+v\left(\frac{1-e^{-\beta t}}{\beta}\right),\\
\dot{X}_{t}(x,v) & = &
v e^{-\beta t},
\end{eqnarray*}
and $\Phi_t(x,v) = [x+vt,v],$ for $\beta=0.$
The asymptotic behaviour of the functions
$\lambda(t),\,\nu(t),\,\mu(t)$ for small $t$ is described (also for $\beta=0$) by
\begin{displaymath}
\begin{array}{r@{\;\;=\;\;}l@{\quad\sim\quad}r@{\quad}r}
\lambda(t) & \alpha t + \sigma \left[\frac{e^{2\beta t} - 4e^{\beta t} + 3}{2\beta^3} + \frac{1}{\beta^2}\:t\right] + \gamma\left[\frac{2}{\beta}t-\frac{2}{\beta^2}(e^{\beta t}-1)\right] & \alpha t, & t\rightarrow 0, \\[0.4cm]
\nu(t) & \sigma \:\frac{e^{2\beta t} - 1}{2\beta} & \sigma t, & t\rightarrow 0, \\[0.4cm]
\mu(t) & \sigma\left( \frac{1- e^{\beta t}}{\beta}\right)^2 + \gamma\:\frac{2(1-e^{\beta t})}{\beta} & - 2\gamma t, & t\rightarrow 0.
\end{array}
\end{displaymath}
And hence: $$
    f(t) \;\::=\;\: 4\lambda(t)\nu(t) - \mu^2(t) \quad\sim\quad 4(\alpha\sigma - \gamma^2)t^2 > 0.
$$

With these preliminaries, the following parabolic reguralization result can be deduced.
\begin{prop}
\label{prop:decay}
For each parameter set $\{\alpha, \beta, \gamma, \sigma\},$ there exist two constants\\
 $B=B(\alpha,\beta,\gamma,\sigma)$ and $T_0=T_0(\alpha,\beta,\gamma,\sigma),$ such that the solution
 of the linear equation (\ref{eq:linEq}) satisfies
\begin{eqnarray}
\label{eq:decayII-v} \Vert \nabla_v w(t)\Vert_X & \leq &
Bt^{-1/2}e^{\kappa t}\Vert w_0\Vert_X,\quad \forall \,0<t\leq T_0,\\
\label{eq:decayII-x} \Vert \nabla_x w(t)\Vert_X & \leq &
Bt^{-1/2}e^{\kappa t}\Vert w_0\Vert_X,\quad \forall \,0<t \leq T_0,
\end{eqnarray}
for all $w_0\in X.$
\end{prop}

The proof is similar to \cite{Ca} and it will be defered to the Appendix.

\begin{remark}
\label{rem:decay}
(a) Observe that the functions $\nabla_x w, \nabla_v w\in \mathcal{C}((0,\infty);X).$
The local boundedness of $\nabla_x w, \nabla_v w$ on any interval
$(\tau, \tau+T_0]$ follows from (\ref{eq:decayI}) and Prop. \ref{prop:decay}.\\
(b) Note that the strategy of the next section will \emph{not} work in the degenerated parabolic
case $\alpha\sigma-\gamma^2=0,$ since the decay rates of Prop. \ref{prop:decay}
would then be $t^{-3/2}$, which is not integrable at $t=0$. Alternative strategies
for this degenerate case were studied in \cite{ALMS00}.
\end{remark}

\subsection{The non-linear equation: local solution}
\label{subsection:nonlinear} Our aim is to solve the following
non-linear initial value problem
\begin{equation}
\label{eq:nonlinear}
w_{t}(t) \;\:=\;\:  \overline{A} w(t)+
\Theta[V[w(t)]]w(t),\quad\forall\, t>0,\qquad w(t=0)\;\:=\;\:w_0\in X,
\end{equation}
where the pseudo-differential operator $\Theta$ is formally defined by
(\ref{eq:1.4}) and the potential $V[w(t)]$ is the (Newton potential) solution of the
Poisson equation
\begin{equation}
\label{eq:Poisson}
-\Delta_{x}V(t,x) \;\:=\;\: n[w(t)](x)\;\:=\;\:\int_{\R^3} w(t,x,v)\,dv,\qquad x\in\R^3,
\end{equation}
for all $t > 0.$ Actually,
if we assume $w(t)\in X$ for all $t\geq 0,$ then the function
$n[w(t)]$ is well-defined for all $t\geq 0$ (cf. Prop. \ref{prop:density}), and the solution
$V[w(t)]$ satisfies the estimates of
Propositions \ref{prop:potentialI}, \ref{prop:potentialII} for all $t\geq 0.$\\

The Propositions \ref{prop:pert} and \ref{prop:decay} motivate the
definition of the Banach space
\begin{eqnarray}
Y_T &:=&\Big\{\,z\in \mathcal{C}([0,T];X) \,|\, \nabla_v z\in \mathcal{C}((0,T];X)\, \mbox{ with } \Vert \nabla_v z(t)\Vert_{X}\leq C t^{-1/2} \;\;\mbox{for}\; t\in (0,T) \Big \},\nonumber
\end{eqnarray}

endowed with the norm
$$ \Vert z \Vert_{Y_T} \;\::=\;\: \sup_{t \in [0,T]}\Vert
z(t)\Vert_{X} + \sup_{t \in [0,T]}\Vert  t^{1/2}\nabla_v z(t)\Vert_{X},$$
for every fixed $0<T < \infty.$ We shall obtain the (local-in-time) well-posedness result
for the problem (\ref{eq:nonlinear}) by introducing a non-linear
iteration in the space $Y_T,$ with an appropriate (small enough) $T.$\\

For a given $w\in Y_T$ we shall now consider the linear Cauchy
problem for the function $z,$
\begin{equation}
\label{eq:linearpert}
z_{t} \;\:=\;\:  \overline{A} z(t)+
\Theta[V[z(t)]]w(t),\quad\forall\, t\in (0,T],\qquad
z(t=0)\;\:=\;\:w_0\in X,
\end{equation}
with $0<T\leq T_0$ and $T_0$ is defined in Prop. \ref{prop:decay}.
According to Prop. \ref{prop:pert} the (time-dependent) operator $\Theta[V[\ldotp]]w(t)$
is, for each fixed $t\in (0,T_0],$ a well-defined, linear and bounded map on $X,$ which we shall consider
as a perturbation of the operator $\overline{A}.$
\begin{lemma}
\label{lemma:nonlinear} For all $w_0\in X$ and $w\in Y_T,$ with
$T\leq T_0,$ the initial value problem
$$
z_{t} \;\:=\;\: \overline{A}z(t)+ \Theta[V[z(t)]]w(t),\quad\forall\, t\in (0,T],\qquad z(t=0)\;\:=\;\:w_0,
$$
has a unique mild solution $z\in\mathcal{C}([0,T];X),$ which satisfies
\begin{equation}
\label{eq:mild}
z(t)\;\:=\;\:\mathrm{e}^{t\overline{A}}w_0+\int_0^t\!\mathrm{e}^{(t-s)
\overline{A}}\Theta[V[z(s)]]w(s)\,ds,\quad \forall \, t \in [0,T].
\end{equation}
Moreover, the solution $z$ belongs to the space $Y_{T}.$
\end{lemma}
\proof. The first assertion follows directly by applying (a trivial extension of)
Thm. 6.1.2 in \cite{Pazy}: \\
For any fixed $w\in Y_T,$ the
function $g(t,.):=\Theta[V[.]]w(t)$ is a bounded linear operator on $X$ for all $t\in (0,T),$
and it satisfies $g\in L^1((0,T);\mathcal{B}(X))\cap C((0,T];\mathcal{B}(X))$ (by Prop. \ref{prop:pert}).
Moreover, by estimates (\ref{eq:estimate}), (\ref{eq:pert1}), the following
inequalities hold
\begin{eqnarray}
\label{eq:stima-1}
\!\!\!\!\!\!\!\!\!\!\!\!\!\!\Vert z(t)\Vert_{X}&\leq&4 \mathrm{e}^{\kappa t}\Vert w_0\Vert_{X} +
4 \int_0^t\!\mathrm{e}^{\kappa(t-s)} C\{\,\Vert w(s)\Vert_{X}+\Vert
\nabla_v w(s)\Vert_{X}\}\Vert z(s)\Vert_{X}\,ds\\
\label{eq:stima0}
&\leq& 4 \mathrm{e}^{\kappa t}\Vert w_0\Vert_{X} + 4C\mathrm{e}^{\kappa
T}\Vert w\Vert_{Y_T}\!\int_0^t\! (1+s^{-1/2})\Vert
z(s)\Vert_{X}\,ds,
\end{eqnarray}
for all $t \in [0,T].$ Then, by Gronwall's Lemma,
\begin{equation}
\label{eq:stima1}
\Vert z(t)\Vert_{X} \;\:\leq\;\: 4e^{\kappa T}\Vert w_0\Vert_{X}
\left[ 1 + 4C\Vert w\Vert_{Y_T} e^{\left(\kappa T + 4Ce^{\kappa T}\Vert w\Vert_{Y_T}(T + 2T^{1/2})\right)}(t+2t^{1/2}) \right],
\end{equation}
for all $t \in [0,T].$
By differentiating equation
(\ref{eq:mild}) in the  $v$-direction, we obtain
\begin{eqnarray}
\label{eq:stima2a}
\nabla_v z(t)\;\:=\;\:\nabla_v e^{t\overline{A}}w_0+\int_0^t\!\nabla_v
e^{(t-s)\overline{A}}g(s,z(s))\,ds,\quad \forall \, t \in [0,T].
\end{eqnarray}
Using the estimates (\ref{eq:decayII-v}), (\ref{eq:pert1}), and
(\ref{eq:stima1}) then yields
\begin{eqnarray}
\label{eq:stima2}
\!\!\!\!\!\!\!\!\!\!\!\!\!\!\!\!\!\Vert \nabla_v z(t)\Vert_{X} &\leq& Bt^{-1/2}\mathrm{e}^{\kappa t}\Vert w_0\Vert_{X}\nonumber\\
&&  +\: B\Vert w\Vert_{Y_T} \!\int_0^t\!(t-s)^{-1/2}\mathrm{e}^{\kappa(t-s)} C\{\,1 + s^{-1/2}\}\,\Vert z(s)\Vert_{X}\,ds \nonumber\\
&\leq& Bt^{-1/2}\mathrm{e}^{\kappa t}\Vert w_0\Vert_{X} + 4BCe^{2\kappa T}\Vert w_0\Vert_X\Vert w\Vert_{Y_T}\!\left[\frac{}{}\pi + 2t^{1/2}\right. \nonumber\\
&& +\left. 4C\Vert w\Vert_{Y_T}e^{\left(\kappa T + 4Ce^{\kappa T}\Vert w\Vert_{Y_T}(T + 2T^{1/2})\right)}\left(4t^{1/2}+\frac{3}{2}\pi t + \frac{4}{3}t^{3/2}\right)\right],
\end{eqnarray}
for all $t \in [0,T].$ The continuity in time of $\nabla_vz$ can be
derived from \fer{eq:stima2a} by using Remark \ref{rem:decay} and the fact that
$g(t,z(t))\in\mathcal{C}((0,T];X).$ Hence, the function $z$ belongs to the space
$Y_{T}.$ \hfill\framebox[0.2cm]\\
\\

We now define the (affine) linear map $M$
on $Y_T$ (for any fixed $0<T\leq T_0$):
$$ w \longmapsto M w:=z,$$
where $z$ is the unique mild solution of
the initial value problem (\ref{eq:linearpert}).
According to Lemma \ref{lemma:nonlinear}, $z\in Y_T.$ Next we shall show that
$M$ is a strict contraction on a closed subset of $Y_{T},$ for $T$ sufficiently
small. This will yield the local-in-time solution of the non-linear equation
(\ref{eq:nonlinear}).

\begin{prop}
\label{prop:M} For any fixed $w_0\in X,$ let $R>\max \{4,B\}\mathrm{e}^{\kappa}\Vert w_0\Vert_{X},$
with the constant $B$ defined in Prop. \ref{prop:decay}.
Then there exists a $\tau:=\tau
(\Vert w_0\Vert_{X}, B)>0$ such that the map $M,$
\begin{equation}
\label{eq:M}
(M w)(t) \;\:=\;\:\mathrm{e}^{t
\overline{A}}w_0+\int_0^t\!\mathrm{e}^{(t-s)
\overline{A}}\Theta[V[M w(s)]]w(s)\,ds,\quad \forall \, t \in[0,\tau],
\end{equation}
is a strict contraction from the ball of radius $R$ of $Y_{\tau}$ into itself.
\end{prop}
\proof.
By (the proof of) Lemma \ref{lemma:nonlinear}, the function $z=M w\in Y_{\tau}$ satisfies (\ref{eq:stima1}).
Under the assumption $\Vert w\Vert_{Y_{\tau}}\leq R,$
this estimate reads
\begin{eqnarray*}
\Vert Mw(t)\Vert_X & \leq & 4e^{\kappa\tau}\Vert w_0\Vert_X
\left[ 1 + 4CRe^{\left(\kappa\tau + 4CRe^{\kappa \tau}(\tau + 2\tau^{1/2})\right)}(t+2t^{1/2}) \right],\quad\forall t\,\in [0,\tau].
\end{eqnarray*}
If we assume
\begin{equation}
\label{eq:C1}
4e^{\kappa\tau}\Vert w_0\Vert_X
\left[ 1 + 4CRe^{\left(\kappa\tau + 4CRe^{\kappa \tau}(\tau +
2\tau^{1/2})\right)}(\tau+2\tau^{1/2}) \right] \;\:\leq\;\: \frac{R}{3},
\end{equation}
then $\Vert Mw(t)\Vert_X \;\leq\; \frac{R}{3}.$ Similar to (\ref{eq:stima2}) we have
\begin{eqnarray*}
\Vert \nabla_v M w(t)\Vert_{X} &\leq& Bt^{-1/2}\mathrm{e}^{\kappa t}\Vert w_0\Vert_{X} + 4BCRe^{2\kappa\tau}\Vert w_0\Vert_X \!\left[\frac{}{}\pi + 2t^{1/2}\right. \nonumber\\
& &+\left. 4CRe^{\left(\kappa \tau + 4CRe^{\kappa \tau}(\tau + 2\tau^{1/2})\right)}\left(4t^{1/2}+\frac{3}{2}\pi t + \frac{4}{3}t^{3/2}\right)\right]\!.
\end{eqnarray*}
If we assume
\begin{eqnarray}
\!\!\!\!\!\!\!\!\!\!\!\!\!\!\!\!\!\!\! B\mathrm{e}^{\kappa\tau}\Vert w_0\Vert_{X}\!\!\! & +& \!\!\!4BCRe^{2\kappa\tau}\Vert w_0\Vert_X \!\left[\frac{}{}\pi\tau^{1/2} + 2\tau \:+\right.\nonumber\\
\label{eq:C2}
\!\!\! &+& \!\!\!\left. 4CRe^{\left(\kappa \tau + 4CRe^{\kappa \tau}(\tau + 2\tau^{1/2})\right)}\left(4\tau + \frac{3}{2}\pi \tau^{3/2} + \frac{4}{3}\tau^2\right)\right] \;\:\leq\;\: \frac{R}{3},
\end{eqnarray}
then
$$t^{1/2}\Vert\nabla_v M w(t)\Vert_{X} \;\:\leq\;\:  \frac{R}{3}, \quad \forall\,t\in[0,\tau].$$

Let us now choose
\begin{eqnarray}
\label{eq:tau}
\tau & := & \min \left\{ 1, \left (\frac{R/3 - 4e^\kappa\Vert w_0\Vert_{X}}{48CR\Vert w_0\Vert_X e^{2\kappa + 12CRe^{\kappa}}}\right)^{2}\!\!,\right.\nonumber\\
 & & \left.\left ( \frac{R/3-B\mathrm{e}^{\kappa}\Vert w_0\Vert_X }{4BCRe^{2\kappa}\Vert w_0\Vert_X \left[ \pi+2  + 4 CRe^{(\kappa + 12CRe^{\kappa})}(\frac{3}{2}\pi + \frac{16}{3})\right]}\right)^2 {}\right\},
\end{eqnarray}
which is positive since $\max \{4,B\}\mathrm{e}^{\kappa}\Vert w_0\Vert_{X}<R.$
Then, the estimates (\ref{eq:C1}) and (\ref{eq:C2})  hold, and hence the
operator $M$ maps the ball of radius $R$ of $Y_\tau$ into itself. \\


To prove contractivity we shall estimate $\Vert M u- Mw\Vert_{Y_\tau}$ for
all $u,w \in Y_\tau$ with $\Vert u\Vert_{Y_\tau},$ $\Vert w\Vert_{Y_\tau}\le R$. Since

\begin{eqnarray*}
M u(t)-M w(t)&=& \int_0^t\!\mathrm{e}^{(t-s) \overline{A}}\Theta[V[(M u-M
w)(s)]]u(s)\,ds\\
& & + \int_0^t\!\mathrm{e}^{(t-s)\overline{A}}\Theta[V[M w(s)]](u-w)(s)\,ds,\quad \forall\, t\in[0,\tau],
\end{eqnarray*}
by analogous estimates,
\begin{eqnarray*}
\Vert M u(t)- Mw(t)\Vert_{X} &\leq& 4 C R\, \mathrm{e}^{\kappa\tau}\left\{\int_0^t\!
(1+s^{-1/2})\,\Vert (M u- M w)(s)\Vert_{X}\,ds \right.\\
& & \left.+\:\Vert u-w\Vert_{Y_{\tau}}\!\int_0^t\!
(1+s^{-1/2})\,ds \right\},
\end{eqnarray*}
and, by applying Gronwall's Lemma:
\begin{eqnarray*}
\label{eq:estimate1}
\Vert M u(t)-M w(t)\Vert_{X} &\leq& 4 C R \mathrm{e}^{\kappa\tau} \left[ t + 2t^{1/2} +\right.\\
& &\!\!\!\!\!\!\!\!\!\!\!\!\!\!\!\!\!\!\!\!\!\!\!\!\!\!\!\!\!\!\!\! + \left.4CR\mathrm{e}^{\left(\kappa\tau + 4CRe^{\kappa\tau}(\tau+2{\tau}^{1/2})\right)}\left( 2t+2t^{3/2} +\frac{1}{2}t^2\right)\right]\Vert u-w\Vert_{Y_{\tau}},\quad \forall \,t \in [0,\tau].
\end{eqnarray*}
By using $0\leq t\leq\tau\leq 1,$ we obtain
\begin{equation}
\label{eq:first}
\Vert M u(t)-M w(t)\Vert_{X} \;\:\leq\;\: 4C R\,\mathrm{e}^{\kappa}\left [3+ 18CR\mathrm{e}^{\left(\kappa + 12CRe^{\kappa}\right)}\right]{\tau}^{1/2}\Vert u-w\Vert_{Y_{\tau}}.
\end{equation}
Similarly,
\begin{eqnarray*}
\Vert \nabla_v M u(t)- \nabla_v Mw(t)\Vert_{X} &\leq& CBR\mathrm{e}^{k\tau}\!\left\{\int_0^t\!(t-s)^{-1/2}(1+s^{-1/2})\,\Vert
(M u- M w)(s)\Vert_{X}\,ds \right. \\
&& +\!\left.\int_0^t\!(t-s)^{-1/2}(1+s^{-1/2})\,ds \Vert
u-w\Vert_{Y_{\tau}}\right\},
\end{eqnarray*}
and, by using estimate (\ref{eq:estimate1}),
\begin{eqnarray*}
\Vert \nabla_v M u(t)- \nabla_v Mw(t)\Vert_{X} & \leq & CBR\mathrm{e}^{\kappa\tau} \left[ 1 + 4C R\,\mathrm{e}^{\kappa}
\left(3+ 18CR\mathrm{e}^{\left(\kappa + 12CRe^{\kappa}\right)}\right)\tau^{1/2}\right]\\
& & \cdot\left(\pi+2t^{1/2}\right) \Vert u-w\Vert_{Y_{\tau}}, \quad \forall \,t \in [0,\tau].
\end{eqnarray*}
Then, by exploiting $0<\tau\leq 1,$
\begin{eqnarray}
\label{eq:second}
\!\!\!\!\!\!\!\!\!\!\! t^{1/2}\Vert \nabla_v M u(t)-\nabla_v M w(t)\Vert_{X} &\leq&
CBR\mathrm{e}^{\kappa} (\pi +2)\left[\frac{}{} 1 + 4C R\,\mathrm{e}^{\kappa}\right.\nonumber\\
&& \left. \!\cdot\left(3+ 18CR\mathrm{e}^{\left(\kappa + 12CRe^{\kappa}\right)}\right)\frac{}{}\!\!\right]\tau^{1/2}\Vert u-w\Vert_{Y_{\tau}}.
\end{eqnarray}
When choosing $\tau>0$ small enough, estimates (\ref{eq:first}), (\ref{eq:second})
imply
$$\Vert M u-M w\Vert_{\mathcal{C}([0,\tau];X)}\;\:\leq\;\: C \Vert u-w\Vert_{\mathcal{C}([0,\tau];X)},$$
for some $C<1,$ and the assertion is proved.  \hfill\framebox[0.2cm]\\
\begin{corollary}
\label{corollary:solution}
There exists a $t_{\mathrm{max}}\leq\infty$
such that the initial value problem (\ref{eq:nonlinear}) has a
unique \emph{mild} solution $w$ in $Y_T,\; \forall
\,T<t_{\mathrm{max}},$ which satisfies
\begin{equation}
\label{eq:mildSol}
w(t) \;\:=\;\: \mathrm{e}^{t\overline{A}}w_0+\int_0^t\!\mathrm{e}^{(t-s)
\overline{A}}\Theta[V[w(s)]]w(s)\,ds,\quad \forall \, t \in[0,T].
\end{equation}
Moreover, if $t_{\mathrm{max}}<\infty,$
then
$$
\lim_{t \nearrow t_{\mathrm{max}}}\Vert w(t)\Vert_{X} \;\:=\;\: \infty.
$$
\end{corollary}

\proof. The solution of the problem is the fixed point of the map
$M$ previously introduced. By Prop. \ref{prop:M} this solution  exists
for a time interval of length $\tau$ (depending only on $\Vert
w_0\Vert_X$) and it belongs to the space $Y_\tau$. Since, in
particular, $w(\tau)\in X,$
the solution can be repeatedly continued up to the maximal time $t_{\mathrm{max}}.$
It will then belong to $Y_T,$ $\forall T<t_{\mathrm{max}}.$\\
If the second assertion of the corollary would not hold, there
would be a sequence of times $t_n\uparrow t_{\mathrm{max}}$ such
that $\Vert w(t_n)\Vert_X\leq C$ for all $n.$ Then, by solving a
problem with the initial value $w(t_n),$ with $t_n$ sufficiently
close to $t_{\mathrm{max}},$ we would extend the solution up to a
certain time $t_n+\tau(\Vert w(t_n)\Vert_X)>t_{\mathrm{max}}.$
This construction would contradict our definition of
$t_{\mathrm{max}}.$\\ The uniqueness of the mild solution follows
by arguments analogous to those in the proof of Thm. 6.1.4 in
\cite{Pazy}.  \hfill\framebox[0.2cm]\\
\begin{remark}
Note that the last statement in the thesis of the Corollary \ref{corollary:solution}
differs from the standard setting (cf. Thm. 6.1.4 in \cite{Pazy}).
For $t_{\mathrm{max}}<\infty$ we conclude the `explosion' of $w(t),\,t\rightarrow t_{max}$ in $X$
and \textnormal{not only} in $Y_t$.
This is due to the parabolic regularization of the problem (\ref{eq:nonlinear}).
\end{remark}

\section{Global-in-time solution, a-priori estimates}
\setcounter{equation}{0}

In this section we shall  exploit dispersive effects of the free transport equation to derive
an a-priori estimate on the electric field. This is the key ingredient for proving the main
result of the paper, the global solution for the WPFP system:

\begin{theorem}
\label{theorem:main} Let $w_0\in X$ satisfy for some $\omega \in
[0,1)$
\renewcommand{\theequation}{{\bf A}}
\begin{equation}
\label{eq:A} \Big\| \!\int\!w_0(x-\vartheta(t)v,v)\,
dv\Big\|_{L^{6/5}(\R^3_x)} \;\:\leq\;\: C_T \vartheta(t)^{-\omega},
\quad \forall\, t\in(0,T], \quad \forall\,T>0,
\end{equation}
\renewcommand{\theequation}{\arabic{section}.\arabic{equation}}
\addtocounter{equation}{-1}

with $\vartheta(t) := \frac{1-e^{-\beta t}}{\beta}$ for $\beta>0,$
and $\vartheta(t) = t$ for $\beta=0.$ Then the WPFP equation
(\ref{eq:nonlinear}) admits a unique {global-in-time} mild
solution $w\in Y_T,\, \forall \,0<T<\infty.$
\end{theorem}

In order to prove that $t_{\mathrm{max}}=\infty$ , we have to show that
$\|w(t)\|_X$ is finite for all $t\ge0$ (cf. Corollary
\ref{corollary:solution}). To this end, we shall derive a-priori estimates for
$\|w(t)\|_2$ and $\||v|^2w(t)\|_2.$ Thus, the proof of Thm.
\ref{theorem:main} will be a consequence of a series of Lemmata, in particular
of Lemma \ref{lemma:L2conservation} and Lemma \ref{lemma:v2L2conservation}. In the sequel, $w(t)$ denotes the
unique mild solution for $0\leq t\leq T,$ for an arbitrary
$0<T<t_{\mathrm{max}}.$
\begin{lemma}
\label{lemma:L2conservation}
For all $w_0\in X$, the mild solution of the WPFP equation (\ref{eq:nonlinear}) satisfies
\begin{equation} \label{eq:L2conservation}
\Vert w(t)\Vert_2^2 \;\:\leq\;\: e^{3\beta t}\Vert w_0\Vert_2^2,\quad\forall \,t\in[0,T].
\end{equation}
\end{lemma}
\proof. Roughly speaking, this follows from the dissipativity of
the operator $\overline{A}-\frac{3\beta}{2}I$ in $L^2(\R^6)$ (cf.
(\ref{l2dissip})) and the skew-symmetry of the pseudo-differential
operator. However, since we are dealing only with the \emph{mild}
solution of the equation, the proof requires an approximation of
$w$ by classical solutions.\\\ Since the solution satisfies $w\in
Y_T,\; \forall \,T<t_{\mathrm{max}},$ Prop. \ref{prop:pert} shows
that the function $f(t):=\Theta[V[w(t)]]w(t),\;
\,t\in(0,t_{\mathrm{max}})$ is well defined and it is in
$\mathcal{C}((0,t_{\mathrm{max}});X)\cap L^1((0,T);X),\, \forall\,
0<T<t_{\mathrm{max}}.$\\ For $0<T<t_{\mathrm{max}}$ fixed, let us
consider the following linear inhomogeneous problem:
\begin{equation}
\label{eq:linapprox} \frac{d}{d t}y(t) \;\:=\;\:
\overline{A}y(t)+f(t),\quad t\in[0,T],\qquad y(t=0)\;\:=\;\:w_0 \in X.
\end{equation}
Its mild solution in $[0,T]$ is the function $w,$ due to the
uniqueness of the mild solution of problem (\ref{eq:nonlinear}).
For this linear problem, we can apply Thm. 4.2.7 of \cite{Pazy}:
The mild solution $w$ is the uniform limit (on $[0,T]$) of
\emph{classical} solutions of problem (\ref{eq:linapprox}). More
precisely, there is a sequence $\{w_0^n\}_{n\in\N}\subset
\mathcal{D}(\overline{A}),\,w_0^n\rightarrow w_0$ in $X,$ and a
sequence
$\{f_n(t)\}\subset\mathcal{C}^1([0,T];X),\,f_n(t)\rightarrow f(t)$
in $L^1((0,T);X).$ And the classical solutions
$y_n\in\mathcal{C}^1([0,T];X)$ of the corresponding problems
\begin{equation}
\label{eq:linapproxn} \frac{d}{d t}y_n(t) \;\:=\;\:
\overline{A}y_n(t)+f_n(t),\quad t\in[0,T],\qquad
y_n(t=0)\;\:=\;\:w_0^n,
\end{equation} converge in $C([0,T];X)$ to the solution $w$ of
problem (\ref{eq:linapprox}).

We shall need these approximating classical solutions $y_n$ in order to justify
the derivation of the a-priori estimate: Multiplying both sides of
(\ref{eq:linapproxn}) by $y_n(t)$ and  integrating yields
$$
\frac{1}{2}\frac{d}{d t}\Vert y_n(t)\Vert_2^2 \;\:\leq\;\: \frac{3\beta}{2}\,
\Vert y_n(t)\Vert_2^2+\int\!\!\!\int\!y_n(t)f_n(t)\,dx\,dv,
$$
since the operator $\overline{A}-\frac{3\beta}{2}$ is dissipative in
$L^2(\R^6)$ (cf. (\ref{l2dissip})). By integrating in $t$ and letting $n\rightarrow \infty,$ we have
$$
\Vert w(t)\Vert_2^2 \;\:\leq\;\: \Vert w_0\Vert_2^2+ 3\beta\!\int_0^t\!\!\Vert
w(s)\Vert_2^2\,ds +2\!\int_0^t\!\!\int\!\!\!\int\!w(s)f(s) \,dx\,dv\,ds, \quad\forall \,t\in[0,T].
$$
The second integral is equal to zero since the pseudo-differential operator $\Theta$
is skew-symmetric. Hence, applying Gronwall's Lemma yields
\begin{equation}
\Vert w(t)\Vert_2^2 \;\:\leq\;\: e^{3\beta t}\Vert w_0\Vert_2^2,\quad\forall \,t\in[0,T].
\end{equation}
\hfill\framebox[0.2cm]\\\\

In order to recover similar estimates for $\||v|^2w(t)\|_2,$ we first need a-priori bounds for the
self-consistent field $E=\nabla V$. To this end, we are going to
exploit dispersive effects of the free streaming operator. We
shall adapt to the Wigner-Poisson and Wigner-Poisson-Fokker-Planck
problems the strategies introduced for the (classical)
Vlasov-Poisson problem (\cite{LiPe91,Pe96}), and  for the
Vlasov-Poisson-Fokker-Planck problem (\cite{Bou93,Bou95,Ca98}).

\subsection{A-priori estimates for the electric field: the Wigner-Poisson\\\hspace*{1.15cm}  case}

To explain the strategy, we first consider the (simpler) WP
problem: Let us assume that ${w}^{\mbox{\tiny{wp}}}$ is a
``regular'' solution of the WP problem (e.g., let
${w}^{\mbox{\tiny{wp}}}(t)\in L_x^2(H^1_v), \nabla_x
V[{w}^{\mbox{\tiny{wp}}}](t)\in\mathcal{C}_B(\R^3),$ uniformly on
bounded $t$-intervals) for which the Duhamel formula holds:
$$
{w}^{\mbox{\tiny{wp}}}(x,v,t)\;\:= \;\:{w}^{\mbox{\tiny{wp}}}_0(x-tv,v) + \int_{0}^{t}\Bigl(\Theta[V[{w}^{\mbox{\tiny{wp}}}]]{w}^{\mbox{\tiny{wp}}}\Bigr)(x-sv,v,t-s)\,ds.
$$
We formally integrate in $v$:
\begin{eqnarray*}
    n[{w}^{\mbox{\tiny{wp}}}](x,t) & = & \int_{\R^3}\!\! {w}^{\mbox{\tiny{wp}}}_0(x-tv,v)\, d v + \int_{0}^{t}\!\!\int_{\R^3}\Bigl(\Theta[V[{w}^{\mbox{\tiny{wp}}}]]{w}^{\mbox{\tiny{wp}}}\Bigr)(x-sv,v,t-s)\,dv\,ds \\[0.15cm]
    &=:&{n}^{\mbox{\tiny{wp}}}_0(x,t)+{n}^{\mbox{\tiny{wp}}}_1(x,t),
\end{eqnarray*}
and split the self-consistent field accordingly:
\begin{eqnarray}
\label{eq:field0WP}
E_0^{\mbox{\tiny{wp}}}(x,t)&:=&\lambda\frac{x}{|x|^3}\ast_{x} n^{\mbox{\tiny{wp}}}_0(x,t)\;\:=\;\: \lambda\frac{x}{|x|^3}\ast_{x}\int w^{\mbox{\tiny{wp}}}_0(x-tv,v)\, d v,\\
\label{eq:field1WP}
E_1^{\mbox{\tiny{wp}}}(x,t)&:=&\lambda\frac{x}{|x|^3}\ast_{x}\int_{0}^{t}\!\!\int\Bigl(\Theta[V[w^{\mbox{\tiny{wp}}}]]w^{\mbox{\tiny{wp}}}\Bigr)(x-sv,v,t-s)\,dv\,ds,
\end{eqnarray}
with $\lambda=\frac{1}{4\pi}$. \\
Then, we can estimate separately the two terms
$E^{\mbox{\tiny{wp}}}_0(t),\,E^{\mbox{\tiny{wp}}}_1(t)$ by
exploting the properties of the convolution kernel $1/|x|$ (cf.
\cite{LiPe91,Pe96} for VP, \cite{Bou93, Bou95, Ca98} for VPFP,
\cite{ALMS00} for WPFP). To this end, we need an appropriate
redefinition of the pseudo-differential operator $\Theta[V]$ in
(\ref{eq:1.4}). It is inspired by the operator $\nabla_x V \cdotp
\nabla_v w$ in the VP equation that can be recovered from
$\Theta[V]w$ in the semiclassical limit (cf. Remark
\ref{remark:semiclassical}).\\

Let us recall that
$
\Theta[V]w(x,v) = {\cal F}_{\eta\rightarrow v}^{-1}\Bigl(i\,\delta V(x,\eta){\cal F}_{v\rightarrow\eta}w(x,\eta)\Bigr).
$
We can rewrite
\begin{equation}
\label{eq:delta}
\delta V(x,\eta) \;\:=\:\! \int\limits_{x-\eta/2}^{x+\eta/2}\! \!\nabla_xV(z)\cdot dz
 \;\:=\:\!  \int\limits_{-1/2}^{1/2} \!\!\eta \cdot \nabla_xV(x-r\eta)\,dr \;\:=\;\: \eta\, \cdotp W(x,\eta)\,,
\end{equation}
with the vector-valued function
$$
W(x,\eta)\;\::=\:\!\int\limits_{-1/2}^{1/2}\!\nabla_xV(x-r\eta)\,dr\,,\quad \forall \,(x,\eta)\in \R^6\,.
$$
Then, we define the vector-valued operator
\begin{equation}
\Gamma[\nabla_x V] u(x,v)\;\::=\;\:\mathcal{F}^{-1}_{\eta\rightarrow v} \Bigl( W(x,\eta)\mathcal{F}_{v\rightarrow \eta}u (x,\eta)\Bigr).
\label{gammadefinition}
\end{equation}
It holds:
\begin{lemma}
\label{lemma:gamma}
Let $\nabla_x V\in \mathcal{C}_B(\R^3)$. Then
\begin{enumerate}
\item $W(x,\eta)\in \mathcal{C}_B(\R^6),$ $\;\;\| W\|_{\infty}\;\leq\; \| \nabla_x V\|_{\infty};$
\item $\Gamma[\nabla_x V]: L^2(\R^6)\rightarrow L^2(\R^6)$ and, for all $u\in L^2(\R^6)\,,$
$$
\|\Gamma[\nabla_x V]u\|_{L^2(\R^6)}\;\:\leq\;\: \| \nabla_x V\|_{\infty}\|u\|_{L^2(\R^6)};
$$
\item $\Gamma[\nabla_x V]: L^2_x(H^1_v)\rightarrow L^2_x(H^1_v)$ and, for all $u\in L^2_x(H^1_v),$
\begin{equation}
\label{eq:gammaestimate}
\|\Gamma[\nabla_x V]u\|_{L^2_x(H^1_v)}\;\:\leq\;\: \| \nabla_x V\|_{\infty}\|u\|_{L^2_x(H^1_v)}.
\end{equation}
\end{enumerate}
\end{lemma}
\proof.
The first and the second assertion are obvious. For (\ref{eq:gammaestimate}) we use
\begin{equation}
\label{eq:divergamma}
\!\!\!\!\!\!\!\!\partial_{v_j}\Gamma[\nabla_x V]u (x,v)\;\,=\;\,i\mathcal{F}^{-1}_{\eta\rightarrow v} \Bigl(\eta_j W(x,\eta)\mathcal{F}_{v\rightarrow \eta}u (x,\eta)\Bigr)\;\,=\;\,\Gamma[\nabla_x V]\partial_{v_j}u; \;\; j=1,2,3.
\end{equation}
\hfill\framebox[0.2cm]\\

\begin{lemma}
\label{lemma:redefinition}
Let $\nabla_x V\in \mathcal{C}_B(\R^3)$ and $u\in L_x^2(H^1_v)$. Then
\begin{equation}
\label{eq:redefinition}
\Theta[V]u(x,v)\;\:=\;\:\diver_v\left(\Gamma[\nabla_x V]u\right)(x,v)
\end{equation}
\end{lemma}
\proof.
By the definition (\ref{eq:delta}) and Lemma \ref{lemma:gamma},
$$
\|\delta V(\ldotp,\eta)\|_{\infty}\;\:\leq\;\: |\eta|\,\| W(\ldotp,\eta)\|_{\infty}\;\:\leq\;\: |\eta|\,\| \nabla_x V\|_{\infty}.
$$
Thus, $\|\Theta[V]u\|_{L^2(\R^6)}\leq \| \nabla_x V\|_{\infty}\|u\|_{L_x^2(H_v^1)};$ the right hand side of equation (\ref{eq:redefinition}) is also well-defined in $L^2(\R^6)$ by estimate (\ref{eq:gammaestimate}). Equality then follows by equation (\ref{eq:delta}) and
$$
i\mathcal{F}^{-1}_{\eta\rightarrow v} \Bigl(\eta \,\cdotp W(x,\eta)\mathcal{F}_{v\rightarrow \eta}u (x,\eta)\Bigr)\;\:=\;\:\sum\limits_{j=1}^{3}\partial_{v_j}\Gamma_j[\nabla_x V]u (x,v)\;\:=\;\:\diver_v\left(\Gamma[\nabla_x V]u\right)(x,v).
$$
\hfill\framebox[0.2cm]\\

\begin{remark}[The semiclassical limit]
\label{remark:semiclassical}
The correctly scaled version of the pseudo-diffe-rential operator with the reduced Planck constant $\hbar = \frac{h}{2\pi}$ reads
$$
\Theta_{\hbar}[V]w(x,v) \;\: = \;\: \frac{i}{(2\pi)^{3/2}}\int_{\R^3} \frac{V(x+\frac{\hbar}{2}\eta) - V(x - \frac{\hbar}{2}\eta)}{\hbar} {\cal F}_{v\rightarrow\eta}
w(x,\eta)e^{iv\cdot\eta}\,d\eta.
$$
Under the assumptions of Lemma \ref{lemma:redefinition}, we thus have
\begin{eqnarray*}
    {\cal F}_{v\rightarrow \eta} (\Theta_{\hbar}[V]w(x,v))  & = & \frac{i}{\hbar}\,\delta V(x,\hbar\eta){\cal F}_{v\rightarrow\eta}w(x,\eta)\\
    & = & i\,W(x,\hbar\eta)\,\cdot\,\eta\,{\cal F}_{v\rightarrow\eta}w(x,\eta).
\end{eqnarray*}
The limit $\hbar\rightarrow 0$ then yields:
$$
i\,W(x,\hbar\eta)\,\cdot\,\eta\,{\cal F}_{v\rightarrow\eta}w(x,\eta) \,\: \longrightarrow \,\: i\nabla_xV(x) \cdot \eta\,{\cal F}_{v\rightarrow\eta}w(x,\eta)
    \,\:=\,\: {\cal F}^{-1}_{\eta\rightarrow v} \Bigl( \nabla_xV(x) \cdot \nabla_v w(x,v)\Bigr);
$$
and hence
$$
    \Theta_{\hbar}[V]w(x,v) \;\: \longrightarrow \;\:  \nabla_xV(x) \cdot \nabla_v w(x,v) \quad\mbox{in} \;\; L^2(\R^6),
$$
which is the non-linear term in the VP equation.
\end{remark}

Using the redefinition (\ref{eq:redefinition}) of the pseudo-differential operator, and under the additional assumptions $w^{\mbox{\tiny{wp}}}\in H^1_x(L^2_v), \,\Delta V[w^{\mbox{\tiny{wp}}}]\in \mathcal{C}_B(\R^3),$ we have for $s\in\R$
\begin{eqnarray}
\label{eq:ThetaInShiftedConvForm}
\Bigl(\Theta[V[w^{\mbox{\tiny{wp}}}]]w^{\mbox{\tiny{wp}}}\Bigr)(x-sv,v) &=& \diver_v \Bigl( \Gamma[\nabla_x V[w^{\mbox{\tiny{wp}}}]] w^{\mbox{\tiny{wp}}}(x-sv,v)\Bigr) \nonumber\\
&& + \:s\:\diver_x \Bigl( \Gamma[\nabla_x V[w^{\mbox{\tiny{wp}}}]]w^{\mbox{\tiny{wp}}}\Bigr)(x-sv,v).
\end{eqnarray}
Thus, the  field $E_1^{\mbox{\tiny{wp}}}$ in (\ref{eq:field1WP}) can be rewritten as ($j=1,2,3$)
\begin{eqnarray}
\label{eq:field1WPbis}
(E_1^{\mbox{\tiny{wp}}})_j(x,t) &\!:=\!& \lambda\frac{x_j}{|x|^3}\ast_{x}\diver_x \int_{0}^{t}\!\!\!s\int \Bigl( \Gamma[\nabla_x V[w^{\mbox{\tiny{wp}}}]] w^{\mbox{\tiny{wp}}}\Bigr)(x-sv,v,t-s)\,dv\,ds\\
&\!=\!&\lambda\sum_{k=1}^3\frac{-3x_jx_k+\delta_{jk}|x|^2}{|x|^5}\ast_{x}\int_{0}^{t}\!\!\!s\int \Bigl( \Gamma_k[\nabla_x V[w^{\mbox{\tiny{wp}}}]] w^{\mbox{\tiny{wp}}}\Bigr)(x-sv,v,t-s)\,dv\,ds.\nonumber
\end{eqnarray}

The following two lemmata are concerned with giving a meaning to the definition (\ref{eq:field1WPbis}) 
of the field $E_1$, independently of the previous derivation.

\begin{lemma}
\label{lemma:technicalEst}
For all $u\in L^2(\R^6)$ and $E\in L^2(\R^3)$ the following estimate holds
\begin{equation}
\label{eq:technicalEst}
\Bigl\|\int_{\R^3_v} \left( \Gamma[E] u\right)(x-sv,v)\,dv\Bigr\|_{L^2(\R^3_x)} \;\:\leq\;\: Cs^{-3/2}\Vert E\Vert_{2}\Vert u\Vert_{2}\,,\quad \forall\, s>0.
\end{equation}

\end{lemma}
\proof.
Since the operator $\Gamma[\ldotp]$ was originally defined for $E\in \mathcal{C}_B(\R^3),$ we shall first derive (\ref{eq:technicalEst}) for $E\in \mathcal{C}_0^{\infty}(\R^3)$ and conclude by a density argument. \\
By the definition (\ref{gammadefinition}) and by several changes of variables, the following chain of equalities holds:
\begin{eqnarray*}
&&\!\!\!\!\!\!\!\!\!\!\!\!\!\!\!\!\!\!\!\!\!\!\!\!\!\!\!\phantom{(2\pi)^3}(\Gamma[E] u)(x,v)
\;\:=\;\: (2\pi)^{3/2}\left[ \mathcal{F}^{-1}_{\eta\rightarrow v} \left(W(x,\eta)\right)\ast_v u\right] (x,v) \\
\!\!\!&\!\!=& \!\!\!\!\int\!\!\!\int\!\!\!\int\limits_{-1/2}^{1/2}\!\! E(x - r\eta)\, e^{i\eta\cdotp z}\,dr\,d\eta \,u(x,v-z) \,dz
\;\:= \;\int\!\!\!\int\!\!\!\int\limits_{-1/2}^{1/2}\!\!\! \frac{1}{|r|^3}E(x - \widetilde{\eta})\, e^{i\tilde{\eta}\cdotp \frac{z}{r}}\,dr\,d\tilde{\eta} \,u(x,v-z) \,dz \\
\!\!\!&\!\!=& \!\!\!\!\int\!\!\!\int\!E(x- \tilde{\eta})\, e^{i\tilde{\eta}\cdotp \tilde{z}}\, d\tilde{\eta}\! \int\limits_{-1/2}^{1/2}\!\!\! \,u(x,v-r\tilde{z}) \,dr\,\,d\tilde{z}
\;\:=\:\int\!\!\!\int\!E(-\hat{\eta})\, e^{i\hat{\eta}\cdotp \tilde{z}}\, d\hat{\eta}\, e^{ix\cdotp \tilde{z}}\!\!\int\limits_{-1/2}^{1/2}\!\!\! \,u(x,v-r\tilde{z}) \,dr\,\,d\tilde{z}\\
\!\!\!&\!\!=& \!\!(2\pi)^{3/2}\int\!\mathcal{F}_{\eta\rightarrow \tilde{z}}E(\tilde{z})\!\int\limits_{-1/2}^{1/2}\!\!\! \,u(x,v-r\tilde{z}) \,dr\, e^{ix\cdotp \tilde{z}}\,d\tilde{z}.
\end{eqnarray*}
Hence
\begin{eqnarray*}
\int \!\left(\Gamma[E] u\right)(x-sv,v)\,dv \!\!&=&\!\!(2\pi)^{\frac{3}{2}}\int\!\mathcal{F}_{\eta\rightarrow \tilde{z}}E(\tilde{z})\left(\int\!\!\!\int_{-1/2}^{1/2}\!\!\! \,u(x-sv,v-r\tilde{z}) \,dr\, e^{-isv\cdotp \tilde{z}}\,dv\!\right)\! e^{ix\cdotp \tilde{z}}\,d\tilde{z}\\
\!\!&=&\!\!\frac{1}{(2\pi s)^3}\int\!\mathcal{F}_{\eta\rightarrow
\tilde{z}}E(\tilde{z})\,\mathcal{F}_{v\rightarrow
\tilde{z}}\!\left(\int_{-1/2}^{1/2}\!\!\!
\,u(x-v,\frac{v}{s}-r\tilde{z}) \,dr\!\right)\! e^{ix\cdotp
\tilde{z}}\,d\tilde{z}\,.
\end{eqnarray*}
Then,
\begin{eqnarray*}
\Big\|\!\int \!\!\left(  \Gamma[E] u\right)(x-sv,v)\,dv\,\Big\|_{L^2(\R^3_x)} &\leq&\!
\frac{\Vert E\Vert_2}{(2\pi s)^3}\Big\|\!\int_{-1/2}^{1/2}\!\!\!\mathcal{F}_{v\rightarrow \tilde{z}}\Bigl(u(x-v,\frac{v}{s}-r\tilde{z})\Bigr)\,dr\Big\|_{L^2(\R^3_{\tilde{z}}\times\R^3_x)}\\
&\leq&\!\frac{\Vert E\Vert_2}{(2\pi s)^3}\!\left(\int_{-1/2}^{1/2}\!\!\big\|\mathcal{F}_{v\rightarrow
\tilde{z}}\Bigl(u(x-v,\frac{v}{s}-r\tilde{z})\Bigr)\big\|_{L^2(\R^3_{\tilde{z}}\times\R^3_x)}^2dr\!\right)^{\!\!\frac{1}{2}}\!\!,
\end{eqnarray*}
by applying  H\"older's inequality first in the $\tilde{z}$ integral and then in the $r$ integral.
Finally, it remains to prove that
$$
\int_{-1/2}^{1/2}\!\Vert\mathcal{F}_{v\rightarrow z}\Bigl(u(x-v,\frac{v}{s}-rz)\Bigr)\Vert_{L^2(\R^{3}_x\times\R^3_z)}^2 dr\;\:=\;\:  s^3\Vert u(x,v)\Vert_{L^2(\R^{3}_{x}\times\R^3_v)}^2.
$$
This is obtained by using repeatedly Plancherel's equality:
\begin{eqnarray*}
&&\!\!\!\!\!\!\!\!\!\!\!\!\!\!\int\limits_{-1/2}^{1/2}\!\!\Vert\mathcal{F}_{v\rightarrow z}\Bigl(u(x-v,\frac{v}{s}-rz)\Bigr)\Vert_{L^2_{x,z}}^2 dr
\:=\:\!\!\int\limits_{-1/2}^{1/2}\!\!\Vert\mathcal{F}_{x\rightarrow \xi}\Bigl[\mathcal{F}_{v\rightarrow z}\Bigl( e^{-iv\xi}u(x,\frac{v}{s}-rz)\Bigr)\Bigr]\Vert_{L^2_{\xi,z}}^2\!dr\\
&\!\!=&\!\!\!\!\!\! \int\limits_{-1/2}^{1/2}\!\!\! \Vert \mathcal{F}_{x\rightarrow \xi}\Bigl(s^3 e^{-is(\xi+z)rz}
\mathcal{F}_{v\rightarrow s(\xi+z)}u(x,v)\Bigr)\Vert_{L^2_{\xi,z}}^2\!dr
\:= \:s^6\!\!\!\int\limits_{-1/2}^{1/2}\!\!\! \Vert \mathcal{F}_{x\rightarrow \xi}\Bigl(
\mathcal{F}_{v\rightarrow s(\xi+z)}u(x,v)\Bigr)\Vert_{L^2_{\xi,z}}^2\!dr \\
&\!\!=&\!\! s^{3}\Vert u(x,v)\Vert_{L^2_{x,v}}^2.
\end{eqnarray*}
\hfill\framebox[0.2cm]\\
\begin{remark}
Observe that the exponent of the variable $s$ recovered in the
Lemma is the same as obtained for the VP case (cf. \cite{Pe96},
e.g.) in the $L^2$-estimate of $\int_{\R^3_v} E u(x-sv,v)\,dv$. In the classical
case, analogous $L^p$-estimates hold in addition. In
the quantum counterpart, instead, the $L^2$-framework is the only possible, since
the estimate had to be derived in Fourier space. Moreover, to derive a more refined version of this basic estimate
(cf. \cite{LiPe91,Bou93}), the non-negativity of the classical
distribution is a crucial ingredient. And this non-negativity does
not hold for Wigner functions.
\end{remark}
The following lemma is an immediate consequence of Lemma
\ref{lemma:technicalEst}.
We shall need the notation
$$
 V_{T,\omega} := \{ E \in C((0,T]; \: L_x^2(\R^3) \:\: \vline
 \:\: \|E\|_{V_{T,\omega}} < \infty \}
$$
with
$$
 \|E\|_{V_{T,\omega}} := \sup_{0 < t \le T} t^\omega
 \|E(t)\|_{L^2}.
$$

\begin{lemma}
\label{lemma:WPL2field} For any fixed $T>0$, let $w\in
\mathcal{C}([0,T];L_{x,v}^2),$ and let $\,w_0$ satisfy for some
$\omega \in [0,1):$
\renewcommand{\theequation}{{\bf B}}
\begin{equation}
\label{eq:condition} \Big\| \!\int\!w_0(x-tv,v)\,
dv\Big\|_{L^{6/5}(\R^3_x)} \;\:\leq\;\: C_T t^{-\omega}, \quad
\forall\, t\in(0,T].
\end{equation}
Then, there exists a unique vector-field $E\in
V_{T,\omega-\frac{1}{2}}$ which satisfies the linear equation
\renewcommand{\theequation}{\arabic{section}.\arabic{equation}}
\begin{equation}
\label{eq:Volterra}
E_j(x,t)\,=\,\lambda\sum_{k=1}^3\frac{-3x_jx_k+\delta_{jk}|x|^2}{|x|^5}\ast_{x}\int_{0}^{t}\!\!\!s\int
\Bigl( \Gamma_k[E_0+E] w\Bigr)(x-sv,v,t-s)\,dv\,ds; \; j=1,2,3
\end{equation}
with $E_0$ defined by $(\lambda = \frac{1}{4\pi})$:
$$
E_0(x,t)\;\::=\;\: \lambda\frac{x}{|x|^3}\ast_{x}\int w_0(x-tv,v)\, d v.
$$
\end{lemma}
\proof. (\ref{eq:Volterra}) has the structure of a Volterra
integral equation of the second kind. Hence, we define the
(affine) map $M: V_{T,\omega-\frac{1}{2}} \rightarrow
V_{T,\omega-\frac{1}{2}}$ by
$$
(M E)_j(x,t)\;\::=\;\:\lambda\sum_{k=1}^3\frac{-3x_jx_k+\delta_{jk}|x|^2}{|x|^5}\ast_{x}\int_{0}^{t}\!\!\!s\int \Bigl( \Gamma_k[E_0+E] w\Bigr)(x-sv,v,t-s)\,dv\,ds.
$$
Applying the generalized Young inequality to the definition of $E_0$ yields
\begin{equation}
\label{eq:WP-field0L2} \Vert E_0(t)\Vert_{2}\;\:\leq\;\: C
\Vert\!\int\!w_0(x-vt,v)\, d v\Vert_{L^{6/5}(\R_x^3)}\,,\quad
\forall\,t\in (0,T]\,.
\end{equation}
Thus, by Lemma \ref{lemma:technicalEst}, the second convolution
factor in (\ref{eq:Volterra}) is well-defined and
$$
\Big\|\!\int_{\R^3_v}\!\!\Bigl( \Gamma_k[E_0+E]
w\Bigr)(x-sv,v,t-s)\,dv\Big\|_{L^2(\R^3_x)} \leq
Cs^{-3/2}\Vert(E_0+E)(t-s)\Vert_{2}\Vert w(t-s)\Vert_{2}\,,\,
\forall\, s\in(0,t].
$$
By classical properties of the convolution with $\frac{1}{|x|}$ (cf. \cite{Stein}) and the Young inequality, we get
\begin{equation}
\label{eq:crucial} \!\!\Vert (M E)_j(t)\Vert_2 \;\leq\; C\int_0^t
\!\frac{1}{\sqrt{s}}\,(\Vert E_0(t-s)\Vert_2+\Vert
E(t-s)\Vert_2)\Vert w(t-s)\Vert_2\,ds\,,\quad \forall\, t\in
(0,T].
\end{equation}
Hence, the map $M$ is well-defined from $V_{T,\omega-\frac{1}{2}}$
into itself and satisfies
\begin{eqnarray*}
\Vert M E (t)\Vert_2 &\leq& C \Bigl(C_T+ \sup_{s\in (0,T]}
s^{\omega-\frac{1}{2}} \Vert E(s)\Vert_2\Bigr)\sup_{s\in
[0,T]}\Vert w(s)\Vert_2\,
\left(t^{1-\omega}+t^{\frac{1}{2}-\omega}\right),\quad \forall
\, t\in(0,T].
\end{eqnarray*}

Since the map is affine, we have (by induction) for all $t\in
(0,T]$
\begin{eqnarray*}
\Vert M^n E(t)- M^n \widetilde{E}(t) \Vert_2 \!\!&\leq&\!\!
\!C\!\sup_{s\in [0,T]}\Vert w(s)\Vert_2 \int_0^t
\!\frac{1}{\sqrt{t-s}}\,
\Vert M^{n-1}E(s)-M^{n-1}\widetilde{E}(s)\Vert_2 \,ds \\
\!\!&\leq&\! \!\!\left(C\!\sup_{s\in [0,T]}\!\Vert
w(s)\Vert_2\right)^{\!n} \!\!C_{n-1}\!\int_0^t
\!\frac{s^{\frac{n}{2}-\omega}}{\sqrt{t-s}}\,ds \sup_{s\in (0,T]}
\!\!\left( s^{\omega - \frac{1}{2}} \Vert E(s)-\widetilde{E}(s)\Vert_2
\right)\!,
\end{eqnarray*}
with
$$
\int_0^t
\!\frac{s^{\frac{n}{2}-\omega}}{\sqrt{t-s}}\,ds \;\:=\;\:
t^{\frac{n+1}{2}-\omega} B\left(\frac{1}{2},
\frac{n+2}{2}-\omega\right)\!, 
$$
$$
\;\;\; C_{n-1} \;\:=\;\:
\prod\limits_{j=1}^{n-1} B\left(\frac{1}{2},
\frac{j}{2}+1-\omega\right) \;\:=\;\: \frac{\pi^{\frac{n-1}{2}}
\Gamma\left( \frac{3}{2}-\omega \right) }{\Gamma\left( \frac{n}{2}
+ 1 - \omega \right)}\,,
$$
where $B$ denotes the Beta function and $\Gamma$ the Gamma
function. Thus, the map $M^n$ is contractive for $n$ large enough
and admits a unique fixed point $E\in V_{T,\omega-\frac{1}{2}}$.
\hfill\framebox[0.2cm]\\
\\

With $E=E_1^{w p}$ the above lemma yields the regularity of the
self-consistent field in the WP equation: It satisfies $\nabla_x
V[w^{\mbox{\tiny{wp}}}]= E_1^{w p}+E^{\mbox{\tiny{wp}}}_0\in
V_{T,\omega - \frac{1}{2}},$ under the assumptions that
$w^{\mbox{\tiny{wp}}}\in \mathcal{C}([0,T];L_{x,v}^2)$ and
$w^{\mbox{\tiny{wp}}}_0$ satisfies (\ref{eq:condition}).
\begin{prop}
\label{prop:WPL2fieldbis} For any fixed $T>0$, let
$w^{\mbox{\tiny{wp}}}\in \mathcal{C}([0,T];L^2_{x,v})$ be a mild
solution of the WP equation with $\Vert
w^{\mbox{\tiny{wp}}}(t)\Vert_2 \;=\; \Vert
w^{\mbox{\tiny{wp}}}_0\Vert_2,$ and with the initial value
$w^{\mbox{\tiny{wp}}}_0$ satisfying condition
(\ref{eq:condition}). Then, the self-consistent field satisfies
the following estimates for all $t\in (0,T]:$
\begin{eqnarray}
\label{eq:WP-field0Lp}
\Vert E^{\mbox{\tiny{wp}}}_0(t)\Vert_{2}&\leq& C \Big\|\!\int\! w^{\mbox{\tiny{wp}}}_0(x-vt,v)\, d v\Big\|_{L^{6/5}(\R_x^3)} \le C C_T t^{-\omega}, \\
\label{eq:WP-field1Lp} \Vert
E^{\mbox{\tiny{wp}}}_1(t)\Vert_{2}&\leq& C\left(\Vert
w^{\mbox{\tiny{wp}}}_0\Vert_2,\mathrm{sup}_{s\in(0,T]}\left\{
s^\omega \Big\|\!\int \! w^{\mbox{\tiny{wp}}}_0(x-sv,v)
dv\Big\|_{L^{6/5}}\right\},T\right)t^{\frac{1}{2}-\omega}.
\end{eqnarray}
Here and in the sequel, the $T$-dependence of the constants $C$ is continuous (on $T\in \R^+$).
\end{prop}
\proof. The first estimate is (\ref{eq:WP-field0L2}) in Lemma
\ref{lemma:WPL2field}. To derive the second one, we exploit eq.
(\ref{eq:crucial}), the conservation of the $L^2$-norm of the
solution and (\ref{eq:WP-field0Lp}):
\begin{eqnarray*}
\Vert E^{\mbox{\tiny{wp}}}_1(t)\Vert_2 &\leq& C\int_0^t \!s^{-1/2}(\Vert  E^{\mbox{\tiny{wp}}}_0(t-s)\Vert_2+\Vert  E^{\mbox{\tiny{wp}}}_1(t-s)\Vert_2)\Vert  w^{\mbox{\tiny{wp}}}(t-s)\Vert_2\,ds\,\\
&\leq& C \Vert w^{\mbox{\tiny{wp}}}_0\Vert_2\,\mathrm{sup}_{s\in(0,T]}\left\{ s^\omega \Big\|\!\int \! w^{\mbox{\tiny{wp}}}_0(x-sv,v) dv\Big\|_{L^{6/5}}\right\}t^{\frac{1}{2}-\omega}\\
&& +\,C \Vert  w^{\mbox{\tiny{wp}}}_0\Vert_2\int_0^t \!(t-s)^{-1/2}\Vert  E^{\mbox{\tiny{wp}}}_1(s)\Vert_2\,ds. \\
\end{eqnarray*}
The thesis follows by Gronwall's Lemma.
\hfill\framebox[0.2cm]\\
$\phantom{x}$\\

We shall now state a simple condition on $w_0$ that implies both
conditions (A), (B). For $w_0 \in L_x^1(L_v^{{6}/{5}})$ a
Strichartz inequality for the free transport equation (cf. Th. 2
in \cite{CaPe96}) reads:
\begin{equation}
 \label{eq:420}
 \left\| \int w_0(x-tv,v) \: dv \right\|_{L^{6/5}(\R^3_x)} \leq
 t^{-\frac{1}{2}}\|w_0\|_{L^1_x(L_v^{6/5})}, \:\:t>0,
\end{equation}
and hence (A) and (B) hold.
\begin{remark}
Let us again compare the a-priori bounds (\ref{eq:WP-field0Lp}),
(\ref{eq:WP-field1Lp}) with their classical counterparts. Using
(\ref{eq:420}) we obtain the same $t^{-\frac{1}{2}}$--singularity
of $\|E^{wp}(t)\|_2$ for the Wigner-Poisson system, as it was
obtained in \cite{CaPe96} for the VP equation. In the latter case,
similar $L^p$-estimates hold for $p$ in a non-trivial interval.
One crucial reason for this difference is the conservation of
$L^p$-norm of the solution: while the WP equation only conserves
the $L^2$-norm, all $L^p$-norms are constant in the VP case.
A second reason is that we cannot exploit any
pseudo-conformal law for the quantum case, since the Wigner functions are
not non-negative (cf. \cite{Pe96}
for the classical case).
\end{remark}
As a by-product we obtain the following result for the self-consistent potential $V$, which follows directly from the splitting $V^{\mbox{\tiny{wp}}}=V^{\mbox{\tiny{wp}}}_0+V^{\mbox{\tiny{wp}}}_1$
\begin{eqnarray}
\label{eq:potential0WP}
V_0^{\mbox{\tiny{wp}}}(x,t)&:=&\lambda\sum_{i=1}^3\frac{x_i}{|x|^3}\ast_{x} {(E_0^{\mbox{\tiny{wp}}})}_i(x,t),\\
\label{eq:potential1WP}
V_1^{\mbox{\tiny{wp}}}(x,t)&:=&\lambda\sum_{i=1}^3\frac{x_i}{|x|^3}\ast_{x}{(E_1^{\mbox{\tiny{wp}}})}_i(x,t).
\end{eqnarray}

\begin{corollary}
\label{cor:WP-potential}
Under the assumptions of Proposition \ref{prop:WPL2fieldbis}, the self-consistent potential $V^{\mbox{\tiny{wp}}}=V^{\mbox{\tiny{wp}}}_0+V^{\mbox{\tiny{wp}}}_1$ satisfies
the following estimates for all $t\in (0,T]:$
\begin{eqnarray}
\label{eq:WP-potential0Lp}
\Vert V^{\mbox{\tiny{wp}}}_0(t)\Vert_{6}&\leq& C C_T t^{-\omega}, \\
\label{eq:WP-potential1Lp} \Vert
V^{\mbox{\tiny{wp}}}_1(t)\Vert_{6}&\leq& C\left(\Vert
w^{\mbox{\tiny{wp}}}_0\Vert_2,\mathrm{sup}_{s\in(0,T]}\left\{
s^\omega \Big\|\!\int \! w^{\mbox{\tiny{wp}}}_0(x-sv,v)
dv\Big\|_{L^{6/5}}\right\},T\right)t^{\frac{1}{2}-\omega}.
\end{eqnarray}
\end{corollary}

\subsection{A-priori estimates for the electric field: the WPFP case}

According to Corollary \ref{corollary:solution}, the mild solution of the WPFP problem satisfies for all $t\in [0,T]$ ($0<T<t_{\mathrm{max}}$)
\begin{eqnarray*}
    w(x,v,t) &=& \!\int\!\!\!\int\!G(t,x,v,x_0,v_0)w_0(x_0,v_0)\, dx_0\,d v_0 \\
    && + \int_{0}^{t}\!\!\int\!\!\!\int\! G(s,x,v,x_0,v_0)(\Theta[V]w)(x_0,v_0,t-s) \,dx_0\,dv_0\,ds
\end{eqnarray*}
with the Green's function $G$ from (\ref{eq:Green}). According to
\cite{SCDM02} we have
$$
    \int_{\R^3}G(t,x,v,x_0,v_0)\,dv \;\:=\;\: R(t)^{-3/2}\mathcal{N}\left(\frac{x-x_0-\vartheta(t) v_0}{\sqrt{R(t)}}\right),
$$
with
\begin{eqnarray}
\mathcal{N}(x) &:=& (2\pi)^{-3/2}\exp\left(-\frac{|x|^2}{2}\right), \\
\label{eq:vartheta}
\vartheta(t) &=& \frac{1-e^{-\beta t}}{\beta} \;\:=\;\: \mathcal{O}(t),\quad  \mbox{for}\;\,t\rightarrow 0, \\
\label{eq:R}
R(t) &:=& 2\alpha t + \sigma \left[\frac{ 4e^{-\beta t}-e^{-2\beta t} +2\beta t - 3}{\beta^3}\right] +
4\gamma\left[\frac{ e^{-\beta t} + \beta t-1}{\beta^2}\right] \;\:=\;\: \mathcal{O}(t),\quad \mbox{for}\;\, t\rightarrow 0.\nonumber\\
&&
\end{eqnarray}

By exploiting the redefinition (\ref{eq:redefinition}) of the pseudo-differential operator,
we obtain the following expression for the density $n[w]$\\[-0.1cm]
$$
 \hspace*{-10.5cm}n[w](x,t) \;\: =\;\: \int_{\R^3}\!w(x,v,t)\, dv
$$
\vspace*{-0.6cm}
\begin{eqnarray*}
&=& \frac{1}{R(t)^{3/2}}\int\!\!\!\int\!\mathcal{N}\left(\frac{x-x_0-\vartheta(t) v_0}{\sqrt{R(t)}}\right)w_0(x_0,v_0)\,dx_0\,dv_0\\
& &  +\!\int_{0}^{t}\!\!\frac{1}{R(s)^{3/2}}\!\!\int\!\!\!\int\! \mathcal{N}\left(\frac{x-x_0-\vartheta(s) v_0}{\sqrt{R(s)}}\right)\diver_{v_0} \left(\Gamma[\nabla_{x_0} V ]  w\right)(x_0,v_0,t-s) \,dx_0\,dv_0\,ds\\
&=& \frac{1}{R(t)^{3/2}}\int\!\!\!\int\!\mathcal{N}\left(\frac{x-x_0}{\sqrt{R(t)}}\right)w_0(x_0-\vartheta(t) v_0,v_0)\,dx_0\,dv_0\\
& & +\!\int_{0}^{t}\!\frac{\vartheta(s)}{R(s)^{2}}\!\!\int\!\!\!\int\! (\nabla_x\mathcal{N})\left(\frac{x-x_0-\vartheta(s) v_0}{\sqrt{R(s)}}\right) \cdotp\left( \Gamma[\nabla_{x_0} V] w\right )(x_0,v_0,t-s) \,dx_0\,dv_0\,ds\\[2mm]
&=& n_0(x,t) + n_1(x,t),
\end{eqnarray*}
where
\begin{eqnarray*}
n_0(x,t)&:=& \frac{1}{R(t)^{3/2}}\,\mathcal{N}\left(\frac{x}{\sqrt{R(t)}}\right)\!\ast_{x}n_0^{\vartheta}(x,t)\,,\quad \mbox{with}\quad n_0^{\vartheta}(x,t):=\int\!\!w_0(x-\vartheta(t) v,v)\,dv, \\
n_1(x,t)&:=& \int_{0}^{t}\!\frac{\vartheta(s)}{R(s)^{3/2}}\,\mathcal{N}\left(\frac{x}{\sqrt{R(s)}}\right)\ast_{x}\diver_x\int\!\Bigl(\Gamma[\nabla_{x} V] w\Bigr )(x-\vartheta(s) v,v,t-s) \,dv\,ds .
\end{eqnarray*}
Correspondingly, we can split the field (with
$\lambda=\frac{1}{4\pi}$):
\begin{eqnarray}
\label{eq:field0}
E_0(x,t)&:=&\lambda\frac{x}{|x|^3}\ast_{x}n_0(x,t) \;\:=\;\: \frac{1}{R(t)^{3/2}}\,\mathcal{N}\left(\frac{x}{\sqrt{R(t)}}\right)\!\ast_{x}\!E_0^{\vartheta}(x,t),\\
\mbox{ with} \!\!\!&& \!\!\! E_0^{\vartheta}(x,t) \;\::=\;\: \lambda\frac{x}{|x|^3}\ast_{x}n^\vartheta_0(x,t), \nonumber\\
\label{eq:field1}
E_1(x,t)&:=&\lambda\frac{x}{|x|^3}\ast_{x}n_1(x,t).
\end{eqnarray}
\begin{remark}
Note that the splitting of the density (and of the electric field)
is the same as in \cite{Bou93,Bou95,Ca98}: in the WPFP case the
two components of the decomposition ($n_0,n_1,$ as well as
$E_0,E_1$) are smoothed versions (in fact, convolutions with a
Gaussian) of those appearing in the WP case (namely
$n^{\mbox{\tiny{wp}}}_0,n^{\mbox{\tiny{wp}}}_1,E^{\mbox{\tiny{wp}}}_0,E^{\mbox{\tiny{wp}}}_1$).
Actually, the density $n_0^{\vartheta}(x,t)$ (which is convoluted
with the Gaussian to give $n_0$) already differs from
$n^{\mbox{\tiny{wp}}}_0(x,t)$ in the WP case because the shift
contains the function $\vartheta$, which is due to friction (and
analogously for $E_0^{\vartheta}(x,t)$ and
$E^{\mbox{\tiny{wp}}}_0(x,t)$).
\end{remark}
{}From Lemma \ref{lemma:technicalEst} we directly get
\begin{equation}
\label{eq:technicalEstbis}
\Big\|\!\int_{\R^3_v} \left( \Gamma[E] u\right)(x-\vartheta(s)v,v, t-s)\,dv\Big\|_{L^2(\R^3_x)} \;\:\leq\;\: C\vartheta(s)^{-3/2}\Vert E(t-s)\Vert_{2}\Vert u(t-s)\Vert_{2}\,,\quad \forall\,t\geq s>0.
\end{equation}

To derive an $L^2$-estimate on the field we shall proceed as in the WP case (Lemma
\ref{lemma:WPL2field}, Proposition \ref{prop:WPL2fieldbis}).

\begin{lemma}
\label{lemma:field1} Let $w$ be the mild solution of the WPFP
equation (\ref{eq:nonlinear}) and let  $w_0\in X$ satisfy
(\ref{eq:A}) for some $\omega \in [0,1)$. For any fixed $T>0$ the
electric field then satisfies $\nabla_x V\in
V_{T,\omega-\frac{1}{2}}$ and the following estimates hold:
\begin{enumerate}
\item for $2\leq p\leq 6,$ $\theta=\frac{3(p-2)}{2p}$
\begin{equation}
\label{eq:field0Lp} \Vert E_0(t)\Vert_{p} \;\:\leq\;\: C(T)\Vert
w_0\Vert_X^\theta \Vert
n_0^{\vartheta}(t)\Vert_{L^{6/5}}^{1-\theta}=\mathcal{O}(t^{-\omega(1-\theta)})
,\quad \forall\,t\in(0,T];
\end{equation}
\item
\begin{equation}
\label{eq:field1L2} \Vert E_1(t)\Vert_{2}\;\:\leq\;\:
C\left(T,\Vert w_0\Vert_2\,, \sup_{s\in(0,T]}\left\{
\vartheta(s)^\omega \Vert n_0^{\vartheta}(s)
\Vert_{L^{6/5}}\right\}\right){t^{\frac{1}{2}-\omega}}\!,\quad
\forall\, t\in (0,T].
\end{equation}
\end{enumerate}
\end{lemma}
\proof.
The estimate for $\Vert E_0(t)\Vert_{p},$ $p\in[2,6]$ is obtained by applying first the generalized Young inequality and then the Young inequality to the expression (\ref{eq:field0})
\begin{eqnarray*}
\Vert E_0(t)\Vert_p&\leq& C \bigg\|\frac{1}{R(t)^{3/2}}\mathcal{N}\left(\frac{x}{\sqrt{R(t)}}\right)\ast_{x}n_0^{\vartheta}(x,t)\bigg\|_q\\
&\leq& C \left\| \frac{1}{R(t)^{3/2}}\mathcal{N}\left(\frac{x}{\sqrt{R(t)}}\right)\right\|_1\Vert n_0^{\vartheta}(x,t)\Vert_q\\
&=& C \Vert n_0^{\vartheta}(t)\Vert_q, \quad\mbox{with}\quad q=
\frac{3p}{p+3}\in [6/5,2].
\end{eqnarray*}

Next we interpolate $n_0^{\vartheta}$ between $L^2$ and $L^{6/5}$,
use (\ref{eq:density1}) and the dissipativity of the operator
$-v\cdotp \nabla_x-\frac{3}{2}$ in $X$ (cf. Lemma
\ref{lemma:dissipA}):
\begin{eqnarray*}
\Vert n_0^{\vartheta}(t)\Vert_q &\leq& C\Vert w_0(x-\vartheta(t)v,v)\Vert_X^{\theta}\Vert n_0^{\vartheta}(t)\Vert_{{6/5}}^{1-\theta}\\
&\leq& C\, e^{\frac{3}{2}\theta\vartheta(t)}\Vert w_0\Vert_X^{\theta}\Vert n_0^{\vartheta}(t)\Vert_{{6/5}}^{1-\theta},
\end{eqnarray*}
with $\theta=\frac{5}{2}-\frac{3}{q}$. Hence
$$
\Vert E_0(t)\Vert_{p} \;\:\leq\;\: C(T)\Vert w_0\Vert_X^\theta
\Vert n_0^{\vartheta}(t)\Vert_{L^{6/5}}^{1-\theta}.
$$
We rewrite the function $E_1(x,t)$ as
\begin{eqnarray}
\label{eq:Volterrabis}
\!\!\!\!\!\!\!\!\!\!\!\!\!\!\!\!(E_1)_j(x,t)
&=&\lambda\sum_{k=1}^3\frac{-3x_jx_k+\delta_{jk}|x|^2}{|x|^5}\ast_{x}\!\int_{0}^{t}\!\frac{\vartheta(s)}{R(s)^{3/2}}\,\mathcal{N}\left(\frac{x}{\sqrt{R(s)}}\right)\ast_{x}F_k(x,t,s)\,ds,
\end{eqnarray}
$$
\mbox{with }\quad F_k(x,t,s)\;\::=\;\:\int\!\left(\Gamma_k[E_0 + E_1] w\right)(x-\vartheta(s) v,v,t-s) \,dv,\nonumber
$$
For estimating it we exploit classical properties of the
convolution with the kernel $\frac{1}{|x|}$ and apply the Young
inequality:
\begin{eqnarray*}
\Vert E_1(t)\Vert_2 & \leq & C\!\int_{0}^{t}\!\!\vartheta(s) \,\bigg\| \frac{1}{R(s)^{3/2}} \mathcal{N}\left(\frac{x}{\sqrt{R(s)}}\right)\ast_{x} F(x,t,s)\bigg\|_2 \,ds\\
& \leq & C\!\int_{0}^{t}\!\!\vartheta(s)\,\bigg\| \frac{1}{R(s)^{3/2}} \mathcal{N}\left(\frac{x}{\sqrt{R(s)}}\right)\bigg\|_1\Vert F(x,t,s) \Vert_2\,ds\\
& \leq & C(T)\|w_0\|_2
\!\int_{0}^{t}\frac{\Vert E_0(t-s)\Vert_2+\Vert E_1(t-s)\Vert_2}{\sqrt{\vartheta(s)}}\,ds,
\end{eqnarray*}
where the last inequality follows from (\ref{eq:technicalEstbis})
and the $L^2$--a-priori estimate on the solution $w$ (cf. Lemma
\ref{lemma:L2conservation}). By applying the estimate
(\ref{eq:field0Lp}) to $\Vert E_0(t)\Vert_2,$ we get
\begin{eqnarray}
\label{eq:E1est} 
\Vert E_1(t)\Vert_2 &\leq&
C(T)\|w_0\|_2
\left(\sup_{t\in(0,T]}\left\{ \vartheta(t)^\omega
\Vert n_0^{\vartheta}(t)\Vert_{L^{6/5}}\right\} 
\int_{0}^{t}\! \vartheta(s)^{-\frac{1}{2}}\vartheta(t-s)^{-\omega}\,ds \right.\nonumber\\
& & +\left.\int_{0}^{t}\!\frac{\Vert E_1(t-s)\Vert_2}{\sqrt{\vartheta(s)}}\,ds\right),
\end{eqnarray}
where the function $\vartheta(s)=\mathcal{O}(s)$ as $ s\to 0.$ Thus the integrals are finite. \\
To establish a solution of (\ref{eq:Volterrabis}) we introduce the fixed point map
\begin{eqnarray*}
(ME)_j(x,t) &:=&\lambda\sum_{k=1}^3\frac{-3x_jx_k+\delta_{jk}|x|^2}{|x|^5}\ast_{x} \\
&& \ast_{x}\!\int_{0}^{t}\!\frac{\vartheta(s)}{R(s)^{3/2}}\,\mathcal{N}\left(\frac{x}{\sqrt{R(s)}}\right)\ast_{x} \int\!\left(\Gamma_k[E_0 + E] w\right)(x-\vartheta(s) v,v,t-s) \,dv\,ds.
\end{eqnarray*}
By using $0<\frac{\vartheta(T)}{T}t \leq \vartheta(t),$
$\forall\,t\in(0,T]$ and (\ref{eq:E1est}), a simple fixed point
argument as in the proof of Lemmma \ref{lemma:WPL2field} with the
contractivity estimate:
$$
\Vert M^n E(t)- M^n \widetilde{E}(t) \Vert_2 \leq 
\left(\!C\sqrt{\frac{T}{\vartheta(T)}}\!\Vert
w_0\Vert_2\!\right)^{\!n}  \!t^{\frac{n+1}{2}-\omega}
 \,\frac{\pi^{\frac{n}{2}}\Gamma\left(\frac{3}{2}-\omega\right)}
 {\Gamma\left(\frac{n+3}{2}-\omega\right)}\, \sup_{s\in (0,T]}\!\!
\left( s^{\omega-\frac{1}{2}}\! \Vert E(s)-\widetilde{E}(s)\Vert_2
\right)
$$
shows that the linear equation (\ref{eq:Volterrabis}) has a unique
solution $E_1\in V_{T,\omega-\frac{1}{2}}.$ Hence $\nabla_x V=E_0
+E_1\in V_{T,\omega-\frac{1}{2}}$ and Gronwall's Lemma then yields
estimate (\ref{eq:field1L2}).
\hfill\framebox[0.2cm]\\\\
\begin{remark}
For the derivation of the a-priori bound on $\Vert E\Vert_2$, we did not use
any moments of $w$ (neither in $x$ nor  $v$), nor pseudo-conformal laws (cf.
\cite{Bou93,Bou95,Pe96,Ca98} for the classical analogue, i.e.\ VPFP). In fact,
the latter are not useful in the quantum case, since the Wigner function
typically also takes negative values. Moreover, the convolution with the
Gaussian did not play a role there; the estimate (\ref{eq:field1L2}) relies
just on the dispersive effect of the free-streaming operator. The parabolic
regularization will be exploited in the ``post-processing'' Proposition 
\ref{prop:field2}.
\end{remark}

The above lemma was the first crucial step towards proving global
existence of the WPFP solution. Next we shall extend this
estimates on the field $E$ to a range of $L^p$-norms:
\begin{prop}
\label{prop:field2} Let $w$ be the mild solution of the WPFP
equation (\ref{eq:nonlinear}) and let  $w_0\in X$ satisfy
(\ref{eq:A}) for some $\omega \in [0,1)$. Then, we have for
any fixed $T>0$ and for all $p\in[2,6)$:
\begin{equation}
\label{eq:fieldLp} \Vert E_1(t)\Vert_{L^p}\;\:\leq\;\:
C\left(T,\Vert w_0\Vert_2,\sup_{s\in(0,T]}\left\{\vartheta(s)^\omega
\Vert
n_0^{\vartheta}(s)\Vert_{L^{6/5}}\right\}\right)t^{\frac{3}{2p}-\frac{1}{4}-\omega}
\!,\quad\forall \,t\in (0,T].
\end{equation}
\end{prop}
\proof.
We shall estimate $E_1(t)$ (cf.~(\ref{eq:Volterrabis})) by using classical properties
of the convolution by the kernel $\frac{1}{|x|}$ and the following 
\begin{equation}
\label{eq:N}
\Bigg\Vert \mathcal{N}\left(\frac{x}{\sqrt{R(s)}}\right)\Bigg\Vert_q \;\:=\;\: CR(s)^{\frac{3}{2q}},\;\forall\, 1 \leq q \leq \infty\,.
\end{equation}
Namely,
\begin{eqnarray*}
\Vert E_1[w](t)\Vert_{L^p(\R^3_x)} \!\! &\leq &\!\!C\!\int_{0}^{t}\!\!\vartheta(s) \,\bigg\| \frac{1}{R(s)^{3/2}} \mathcal{N}\left(\frac{x}{\sqrt{R(s)}}\right)\ast_{x} F(x,t,s)\bigg\|_{L^p(\R^3_x)} \,ds\\
&\leq &C\!\int_0^t\!\!\frac{R(s)^{\frac{3}{2q}-\frac{3}{2}}}{\sqrt{\vartheta(s)}}\Bigl(\Vert E_0(t-s)\Vert_2+\Vert E_1[w](t-s)\Vert_2\Bigr)\Vert w(t-s)\Vert_2\,ds,
\end{eqnarray*}
where we used the Young inequality with $1+1/p=1/q+1/2$ (thus, $p\geq 2$) and, for the $L^2$-norm, the Lemma \ref{lemma:technicalEst}.
Then, by applying Lemma \ref{lemma:field1} we get 
\begin{eqnarray*}
\Vert E_1[w](t)\Vert_{L^p(\R^3_x)} \!\! &\leq &\!\! C\!\left(T,\sup_{s\in(0,T]}\left\{\vartheta(s)^{\omega} \Vert n_0^{\vartheta}(s)\Vert_{L^{6/5}}\right\},\Vert w_0\Vert_{L^2(\R^6)}\right).\\
&&\cdotp\!\int_0^t\!\!
\frac{R(s)^{\frac{3}{2q}-\frac{3}{2}}}{\sqrt{\vartheta(s)}}\left(\vartheta(t-s)^{-\omega}+(t-s)^{\frac12-\omega}\right)ds\\
\end{eqnarray*}
Since $\vartheta(t)=\mathcal{O}(t)\,, R(t)=\mathcal{O}(t)$ for $t \rightarrow 0$
(cf.~(\ref{eq:vartheta}), (\ref{eq:R})), the last integral is finite for
all $t>0$ and for $3/(2q)-2>-1\, \Leftrightarrow \, 3/(2p)-5/4>-1\, \Leftrightarrow \,p<6$. In fact the integral is
$\mathcal{O}(t^{\frac{3}{2p}-\frac{1}{4}-\omega})$.
\hfill\framebox[0.2cm]\\

\begin{remark}
Proposition \ref{prop:field2} provides a non-trivial
interval of $L^p$-estimates for the electric field in the
WPFP case. This is due to the regularizing effect
of the FP term. We remark that the corresponding Gaussian
is ``better behaved'' than the classical one, since the
quantum FP operator is uniformly elliptic in both $x$ and
$v$ variables. On the other hand, exactly as in the WP
case, the range of $L^p$-estimates for the WPFP equation
is smaller in comparison to the counterpart VPFP and that
depends again on the non-negativity of the classical
distribution function.
\end{remark}
As a further result, we obtain an a-posteriori information on the self-consistent potential $V$, which follows directly from the a-priori estimates on the field. 
Accordingly, we split the potential as $V=V_0+V_1$, with
\begin{eqnarray}
\label{eq:potential0WPFP}
V_0(x,t)&:=&\lambda\sum_{i=1}^3\frac{x_i}{|x|^3}\ast_{x} {(E_0)}_i(x,t),\\
\label{eq:potential1WPFP}
V_1(x,t)&:=&\lambda\sum_{i=1}^3\frac{x_i}{|x|^3}\ast_{x}{(E_1)}_i(x,t).
\end{eqnarray}

\begin{corollary}
\label{cor:WPFP-potential}
Under the assumptions of Proposition \ref{prop:field2}, the self-consistent potential $V(t)=V_0(t)+V_1(t)$ belongs to $L^p(\R^3_x)$ with $6\leq p \leq \infty,$ and satisfies for all $t\in (0,T]:$  
\begin{eqnarray}
\label{eq:WPFP-potential0Lp}
\!\!\!\!\!\!\!\!\!\!\Vert V_0(t)\Vert_{p} &\leq &C(T)\Vert
w_0\Vert_X^\theta \Vert
n_0^{\vartheta}(t)\Vert_{L^{6/5}}^{1-\theta}=\mathcal{O}(t^{-\omega(1-\theta)})
\quad \mbox{with }\;\theta=\frac12 -\frac{3}{p},\\
\label{eq:WPFP-potential1Lp} 
\!\!\!\!\!\!\!\!\!\!\Vert
V_1(t)\Vert_{p}&\leq&
C\left(T,\Vert w_0\Vert_2\,, \sup_{s\in(0,T]}\left\{
\vartheta(s)^\omega \Vert n_0^{\vartheta}(s)
\Vert_{L^{6/5}}\right\}\right){t^{\frac{3}{2p}+\frac{1}{4}-\omega}}\!.
\end{eqnarray}

\end{corollary}

\subsection{A-priori estimates for the weighted $L^2$-norms}

A first consequence of the a-priori estimates for the electric field is the following
\begin{lemma}
\label{lemma:vL2conservation}
For all $w_0\in X$ such that (\ref{eq:A}) holds for some $\omega\in[0,1)$,
the mild solution of the WPFP equation (\ref{eq:nonlinear}) satisfies
\begin{equation}
\label{eq:vL2conservation}
\Vert vw(t)\Vert_2^2 \;\:\leq\;\: C\!\left(T,\Vert w_0\Vert_X,
\sup_{s\in(0,T]}\left\{ \vartheta(s)^\omega \Vert
n_0^{\vartheta}(s)\Vert_{L^{6/5}}\right\}\right),\quad\forall\, t\in[0,T].
\end{equation}
\end{lemma}

\proof.
In order to justify the derivation of this a-priori estimate we need again the approximating classical solutions $y_n$ introduced in the proof of Lemma \ref{lemma:L2conservation}. Mutiplying both sides of 
(\ref{eq:linapproxn}) 
by $v_i^2y_n(t)$ 
and integrating yields 
$$
\frac{1}{2}\frac{d}{dt}\Vert v_i y_n(t)\Vert_2^2 \;\:=\;\:  \!\int\!\!\!\int\! v_{i}^2y_n(t)\overline{A}y_n(t)\,dx\,dv +\int\!\!\!\int\!v_i^2y_n(t)f_n(t)\,dx\,dv.
$$
By analogous calculations as in the proof of Lemma \ref{lemma:dissipA} (cf. also \fer{eq:vi4uAu}) we get,
$$
 \int\!\!\!\int\! |v|^2y_n(t)\overline{A}y_n(t)\,dx\,dv\;\:\leq\;\: 3\sigma \Vert y_n(t)\Vert_2^2 + \frac{\beta}{2}\Vert v y_n(t)\Vert_2^2,
$$
and hence
$$
\frac{1}{2}\frac{d}{d t}\Vert v y_n(t)\Vert_2^2\;\:\leq\;\: 3\sigma \Vert y_n(t)\Vert_2^2 + \frac{\beta}{2}\,\Vert v y_n(t)\Vert_2^2+ \int\!\!\!\int\!|v|^2y_n(t)f_n(t) \,dx\,dv,\quad\forall\,t\in[0,T].
$$
By integrating in $t$, letting $n\rightarrow \infty,$ and using (\ref{eq:L2conservation}), we have
\begin{eqnarray*}
\Vert v w(t)\Vert_2^2 &\leq&\Vert v w_0\Vert_2^2+\frac{2\sigma}{\beta}(e^{3\beta t}-1)
\Vert w_0\Vert_2^2+ \beta\!\int_0^t\!\!\Vert v w(s)\Vert_2^2\,ds \\
&& + \,2\!\int_0^t\!\!\int\!\!\!\int\!|v|^2w(s)f(s)\,dx\,dv\,ds, \quad\forall\,t\in[0,T].
\end{eqnarray*}
Using again the skew-symmetry of the pseudo-differential operator and the H\"older inequality yields
\begin{eqnarray*}
\int_0^t\!\!\int\!\!\!\int\!v_iw(s)v_if(s)\,dx\,dv\,ds &=&
\frac{1}{2}\int_0^t\!\!\int\!\!\!\int\!v_iw(s)\Omega[\partial_iV[w(s)]]w(s)\,dx\,dv\,ds\\
&\leq& \frac{1}{2}\int_0^t\!\!\Vert v_i w(s)\Vert_2\,\Vert\Omega[\partial_iV[w(s)]]w(s)\Vert_2\,ds,
\end{eqnarray*}
with the operator $\Omega$ defined in (\ref{eq:defOmega}). Estimating as in (\ref{eq:deltaplus}) and using
  the Sobolev inequality we obtain for $t\in[0,T]$: 
\begin{eqnarray}
\label{eq:omegaEst}
    \Vert\Omega[\partial_iV[w(t)]]w(t)\Vert_{L^2(\R^3_x\times\R^3_v)} &\leq & C\Vert \partial_iV[w(t)]]\hat{w}(t)\Vert_{L^2(\R^3_x\times\R^3_\eta)} \nonumber\\
    & \leq & C\Vert \partial_iV[w(t)]] \Vert_3 \Vert\hat{w}(t)\Vert_{L^2(\R^3_x; L^6(\R^3_\eta))} \nonumber\\
    & \leq & C\Vert \partial_iV[w(t)]] \Vert_3 \Vert \nabla_\eta\hat{w}(t)\Vert_{L^2(\R^3_x\times \R^3_\eta)} \nonumber\\
    & \leq & C\Vert \partial_iV[w(t)]] \Vert_3 \Vert vw(t)\Vert_2,
\end{eqnarray}
where $\hat{w}(x,\eta,t) := \mathcal{F}_{v\rightarrow \eta}(w(x,v,t)).$
Finally, using Prop. \ref{prop:field2} (estimate (\ref{eq:fieldLp}) with $p=3$) yields
\begin{eqnarray}
\label{eq:v-gronwall}
\Vert v w(t)\Vert_2^2 &\leq& C(T)\left(\Vert v w_0\Vert_2^2+\Vert w_0\Vert_2^2\right) 
+ C\left(T,\Vert w_0\Vert_X,\sup_{s\in(0,T]}\left\{ \vartheta(s)^\omega \Vert
n_0^{\vartheta}(s)\Vert_{L^{6/5}}\right\}\right)\nonumber\\
& & \cdot \!\int_0^t\!\!\left(s^{-\frac{\omega}{2}}+s^{-\omega+\frac{1}{4}} + \beta\right)
\Vert v w(s)\Vert_2^2\,ds,\quad t\in[0,T],
\end{eqnarray}
and the Gronwall Lemma gives the result.
\hfill\framebox[0.2cm]\\\\

With this result we can proceed to derive the a-priori estimate for $\||v|^2w(t)\|_2.$

\begin{lemma}
\label{lemma:v2L2conservation}
For all $w_0\in X$ such that (\ref{eq:A}) holds for some $\omega \in[0,1)$,
the mild solution of the WPFP equation (\ref{eq:nonlinear}) satisfies
\begin{equation} \label{eq:v2L2conservation}
\Vert |v|^2w(t)\Vert_2^2 \;\:\leq\;\: C\!\left(T,\Vert w_0\Vert_X,\sup_{s\in(0,T]}\left\{ \vartheta(s)^\omega \Vert
n_0^{\vartheta}(s)\Vert_{L^{6/5}}\right\}\right),\quad\forall\, t\in[0,T].
\end{equation}
\end{lemma}

\proof.
In order to control the term $\Vert |v|^2w(t)\Vert_2 $, we shall use the same strategy as in the Lemmata \ref{lemma:L2conservation} and \ref{lemma:vL2conservation}. Multiplying both sides of (\ref{eq:linapproxn}) 
by $v_i^4 y_n(t)$ and integrating we get by using (\ref{eq:vi4uAu}) and repeating the same limit procedure as in the previous lemma:
\begin{eqnarray*}
\frac{1}{2}\frac{d}{d t}\sum_{i=1}^{3} \int\!\!\!\int\! v_i^4w(t)^2\,dx\,dv &\leq& 9\sigma \Vert w(t)\Vert_2^2 + \left(3\sigma-\frac{1}{2}\beta\right)\sum_{i=1}^{3}\int\!\!\!\int\! v_i^4 w(t)^2\,dx\,dv\\
&& +\sum_{i=1}^{3}\int\!\!\!\int\!v_i^4 w(t) f(t)\,dx\,dv,\quad\forall\, t\in[0,T].
\end{eqnarray*}
By integrating in $t,$ using $C_1|v|^4 \leq \sum v_i^4 \leq C_2|v|^4$ and (\ref{eq:L2conservation}), we have
\begin{eqnarray}
\label{eq:v2w-est}
\Vert |v|^2w(t)\Vert_2^2 &\leq& C\Biggl(\Vert |v|^2w_0\Vert_2^2 + \frac{6\sigma}{\beta}(e^{3\beta t}-1)\Vert w_0\Vert_2^2 + \left(6\sigma-\beta\right)\!\int_0^t\!\Vert |v|^2w(s)\Vert_2^2\,ds \Biggr.\nonumber\\
&& \left.+2\sum_{i=1}^{3}\int_0^t\!\!\int\!\!\!\int\!v_i^4 w(s) f(s)\,dx\,dv\,ds \right),\quad\forall\, t\in[0,T].
\end{eqnarray}
Using again the skew-symmetry of the pseudo-differential operator $\Theta,$ the equation \fer{eq:2.9} and the H\"older inequality, we have
\begin{eqnarray}
\label{eq:weightedfreeterm}
\int_0^t\!\!\int\!\!\!\int\!v_i^2 w(s) v_i^2 f(s) dx\,dv\,ds &\leq& \frac{1}{4}\!\int_0^t\!\!\Vert v_i^2w(s)\Vert_2\,\Vert\Theta[\partial_i^2V[w(s)]]w(s)\Vert_2\,ds \nonumber\\
& &+\int_0^t\!\!\Vert v_i^2
w(s)\Vert_2\,\Vert\Omega[\partial_iV[w(s)]]v_iw(s)\Vert_2\,ds.
\end{eqnarray}
Since $\hat{w}(x,\ldotp,t)\in H^2(\R_{\eta}^3),$ the Gagliardo-Nirenberg inequality yields for $t\in[0,T]$
\begin{eqnarray}
\label{eq:Linftyw-hat}
\Vert \hat{w}(x,\ldotp,t)\Vert_{L^\infty(\R_{\eta}^3)} &\leq& C\Vert\hat{w}(x,\ldotp,t)\Vert_{L^6(\R^3_\eta)}^{1/2}\Vert \widehat{|v|^2w}(x,\ldotp,t)\Vert_{L^2(\R^3_\eta)}^{1/2} \nonumber\\
&\leq & C\Vert\widehat{v w}(x,\ldotp,t)\Vert_{L^2(\R^3_\eta)}^{1/2}\Vert \widehat{|v|^2w}(x,\ldotp,t)\Vert_{L^2(\R^3_\eta)}^{1/2}.
\end{eqnarray}
Using
$$
\left \Vert \Delta V[w(t)]\right\Vert_2 \;\:=\;\: \left\Vert n[w(t)]\right\Vert_2 \;\:=\;\: C\Vert\hat{w}(\ldotp,\eta=0,t)\Vert_{L^2(\R^3_x)} \;\:\leq\;\: C\!\left(\int\!\!\Vert \hat{w}(x,\ldotp,t)\Vert_{L^\infty(\R^3_\eta)}^2\,dx\right)^{1/2}\!\!,
$$
(\ref{eq:Linftyw-hat}), the H\"older inequality, and (\ref{eq:vL2conservation}) we can estimate:
\begin{eqnarray}
\label{eq:thetaest}
\Vert\Theta[\partial_i^2V[w(t)]]w(t)\Vert_2 &\leq& C\Vert \Delta V[w(t)]\Vert_2 \!\left(\int\!\!\Vert \hat{w}(x,\ldotp,t)\Vert_{L^\infty(\R^3_\eta)}^2\,dx\right)^{1/2}\nonumber\\
&\leq&C \!\int\!\!\Vert\hat{w}(x,\ldotp,t)\Vert_{\infty}^2\,dx \nonumber\\
&\leq& C \!\int\!\! \Vert\widehat{v w}(x,\ldotp,t)\Vert_{L^2(\R^3_\eta)}\;
\Vert \widehat{|v|^2w}(x,\ldotp,t)\Vert_{L^2(\R^3_\eta)}\,dx \nonumber\\
&\leq& 
C\!\left(T,\Vert w_0\Vert_X,\sup_{s\in(0,T]}\left\{ \vartheta(s)^\omega \Vert n_0^{\vartheta}(s)\Vert_{L^{6/5}}\right\}\right)\Vert |v|^2 w(t)\Vert_2.\quad
\end{eqnarray}
For the second term of the r.h.s. of (\ref{eq:weightedfreeterm}) we proceed as in (\ref{eq:omegaEst}) and use the estimate (\ref{eq:fieldLp}):
\begin{eqnarray*}
\hspace*{-4.5cm}\Vert\Omega[\partial_iV[w(t)]]v_iw(t)\Vert_2 &\leq& C\Vert \partial_{i} V[w(t)]\Vert_{3} \Vert \widehat{v_iw}(t)\Vert_{L^2(\R^3_x; L^6(\R^3_\eta))} 
\end{eqnarray*}
\begin{equation}
\label{eq:omegaest}
\leq \;\: C\!\left(T,\Vert w_0\Vert_X,\sup_{s\in(0,T]}\left\{ \vartheta(s)^\omega 
\Vert n_0^{\vartheta}(s)\Vert_{L^{6/5}}\right\}\right)
\left(t^{-\frac{\omega}{2}}+t^{-\omega+\frac{1}{4}}\right)\Vert{|v|^2 w}(t)\Vert_2, 
\quad\forall\, t\in[0,T].
\end{equation}
Analogously to (\ref{eq:v-gronwall}), combining the estimates (\ref{eq:weightedfreeterm}), (\ref{eq:thetaest}) and (\ref{eq:omegaest}) the Gronwall Lemma gives the assertion.
\hfill\framebox[0.2cm]\\\\

{\sc Proof of Theorem \ref{theorem:main}.}\\
The Lemmata \ref{lemma:L2conservation} and \ref{lemma:v2L2conservation} show that
$$\Vert w(t)\Vert_X \;\:\leq\;\: 
C\!\left(T,\Vert w_0\Vert_X,\sup_{s\in(0,T]}\left\{ \vartheta(s)^\omega \Vert
n_0^{\vartheta}(s)\Vert_{L^{6/5}}\right\}\right),\quad\forall\, t\in[0,T], \quad \forall\, 0<T<t_{\mathrm{max}},
$$
with $C$ being continuous in $T\in[0,t_{\mathrm{max}}].$
Then, Corollary \ref{corollary:solution} shows that the mild solution $w$ exists on $[0,\infty)$.
\hfill\framebox[0.2cm]\\\\

\subsection{Regularity}
\label{subsec:regularity}

The following result concerns the smoothness of the solution of WPFP, the macroscopic density and the force field, for positive times. 

\begin{corollary}
\label{corollary:regulOfsolution}
Under the  assumptions of Theorem
\ref{theorem:main},
the mild solution of the WPFP equation
(\ref{eq:nonlinear}) satisfies
$$
w\in \mathcal{C}((0,\infty);\mathcal{C}^\infty_{\mathcal{B}}(\R^6)),
$$

$$
n(t), E(t), V(t)\in \mathcal{C}((0,\infty);\mathcal{C}^\infty_{\mathcal{B}}(\R^3)).
$$
\end{corollary}
\proof. Obviously, $w(t)\in \mathcal{C}(\R^6)$ $\forall t>0,$ because of
the Green's function representation in (\ref{eq:mildSol}), \fer{eq:solution}.
If we differentiate equation (\ref{eq:nonlinear}) with respect to $x_i$ and,
resp., $v_i,$ we obtain the following linear, inhomogeneous problems for
any fixed $t_1>0$.
\begin{eqnarray*}
z_{t}(t) &=&  \overline{A} z(t)+ \Theta[V[z(t)]]w(t) +  \Theta[V[w(t)]]z(t),\quad\forall\, t>t_1,\quad
z(t_1)=\partial_{x_i}w(t_1) \in X,\\
y_{t}(t) &=&  \overline{A} y(t)+ \beta y(t) -
\partial_{x_i}w(t) + \Theta[V[w(t)]]y(t),\quad\forall\,
t>t_1,\quad y(t_1)=\partial_{v_i}w(t_1) \in X.
\end{eqnarray*}
By arguments analogous to Lemma \ref{lemma:nonlinear}, there exists a unique mild solution
\begin{equation}
\label{eq:derivatives}
z=\partial_{x_i}w\in \mathcal{C}([t_1,\infty);H^1(\R^6;(1 + |v|^2)^2\,dx\,dv)).
\end{equation}
By an induction procedure, the derivatives $\nabla^\alpha_x\nabla^\beta_vw,$
for $\alpha,\beta\in\N^3,$ $|\alpha|+|\beta|=m>1$ are also mild solutions of
similar problems with additional well-defined inhomogeneities and with initial
times $0<t_1<t_2< ... <t_m$. This yields  $\nabla^\alpha_x\nabla^\beta_vw\in
\mathcal{C}([t_m,\infty);H^1(\R^6;(1 + |v|^2)^2\,dx\,dv)),$ and thus
$\nabla^\alpha_x\nabla^\beta_vw\in \mathcal{C}((0,\infty);X).$
Hence, the statement about smoothness of the density and the electric field is
straightforward from Propositions \ref{prop:density} and \ref{prop:potentialI}
and Sobolev embeddings.\qed

\section{Appendix}
\setcounter{equation}{0}
\subsection*{Proof of Lemma \ref{lemma:dissipA}}
For $u\in D(A)$ we have
\begin{equation}
\label{innpro}
<Au,u>_{\tilde{X}} \;\:=\;\: <Au,u>_{ L^{2}(\R^6)} +  \sum_{i=1}^{3}\INT v_{i}^4uAu,
\end{equation}
where $\int\!\!\int \!f$ denotes the integral
$\int_{\R^3}\!\int_{\R^3}\!f(x,v)\,dv\,dx\,,$ and the norm $\Vert \ldotp\Vert_{\tilde{X}}$ is defined by (\ref{eq:X-tilde}). Using integrations by parts
we shall calculate the three terms on the right hand side
separately.
\begin{eqnarray*}
 <Au,u>_{ L^{2}(\R^6)} & = & \sum_{i=1}^{3} \left( -\INT v_{i}u_{x_i}u  + \beta\INT (v_{i}u)_{v_i}u \right.\\
& & \left. +\sigma\INT u_{v_i v_i}u + 2\gamma\INT u_{x_i v_i}u + \alpha\INT u_{x_i x_i}u\right) \\
& \leq &  \sum_{i=1}^{3} \left[ 3\beta\INT u^2 + \beta\INT v_i u_{v_i}u - \sigma\INT u_{v_i}^{2} \right.\\
& & \left. +\gamma\left(\epsilon \INT u_{x_i}^{2}
+ \frac{1}{\epsilon}\INT u_{v_i}^{2}\right) - \alpha\INT u_{x_i}^{2}\right] \\
& = & \frac{3}{2}\beta \Vert u\Vert_{2}^{2} + \left(\frac{\gamma}{\epsilon} - \sigma \right)\Vert \nabla_vu \Vert_{2}^{2} + \left( \epsilon\gamma -\alpha\right)\Vert \nabla_xu\Vert_{2}^{2}.
\end{eqnarray*}
With $\epsilon = \frac{\gamma}{\sigma}$ we obtain
\begin{equation}\label{l2dissip}
 <Au,u>_{ L^{2}(\R^6)}  \;\:\leq\;\:  \frac{3}{2}\beta \Vert
 u\Vert_{2}^{2}.
\end{equation}

Next we estimate the second term of \fer{innpro}:
\begin{eqnarray}
\label{eq:vi4uAu}
 \sum_{i=1}^{3} \INT v_{i}^4uAu & = & \sum_{i,j=1}^{3} \left( -\INT v^{4}_{i}v_{j}u_{x_j}u + \beta\INT  v^{4}_{i}(v_{j}u)_{v_j}u \right.\nonumber\\
& & \left. +\sigma\INT  v^{4}_{i}u_{v_j v_j}u + 2\gamma\INT  v^{4}_{i}u_{x_j v_j}u + \alpha\INT  v^{4}_{i}u_{x_j x_j}u\right) \nonumber\\
& \leq &  \sum_{i,j=1}^{3} \left[ \beta\INT  v^{4}_{i}u^2 + \beta\INT v^{4}_{i}v_j u_{v_j}u - \sigma\INT v^{4}_{i}u_{v_j}^{2}  \right.\nonumber\\
& & \left. - \frac{4}{3}\sigma\INT v^{3}_{i}u_{v_i}u + \gamma\left(\epsilon \INT v^{4}_{i}u_{x_j}^{2} + \frac{1}{\epsilon}\INT  v^{4}_{i}u_{v_j}^{2}\right) - \alpha\INT v^{4}_{i}u_{x_j}^{2} \right]\nonumber\\
& \leq & \sum_{i=1}^{3} \left( -\frac{1}{2}\beta\INT v^{4}_{i}u^2 + 6\sigma\INT v_{i}^{2}u^2\right)\nonumber\\
& \leq &  9\sigma\Vert u\Vert_{2}^{2} + \left( -\frac{1}{2}\beta + 3\sigma\right)\sum_{i=1}^{3}\INT v^{4}_{i}u^2,
\end{eqnarray}
by choosing $\epsilon = \frac{\gamma}{\sigma}$ and by an interpolation.

Collecting the two estimates yields
\begin{eqnarray*}
 <Au,u>_{\tilde{X}} & \leq & \left(\frac{3}{2}\beta + 9\sigma\right) \Vert u\Vert_{2}^{2} 
 + 3\sigma\sum_{i=1}^{3}\INT v^{4}_{i}u^2  \;\:\leq\;\: \left( \frac{3}{2}\beta + 9\sigma \right)\Vert u\Vert_{\tilde{X}}^{2}.
\end{eqnarray*}
Thus, the operator $A -  \kappa I$ is
dissipative.\hfill\framebox[0.2cm]\\

\subsection*{Proof of Lemma \ref{lemma:density}}

To prove the assertion we shall construct for each $f\in
D(\overline{P})\subset L^2(\R^6)$ a sequence $\{f_n\}\subset D(P)$
such that $f_n\to f$ in the graph norm \\
$\|f\|_P=\|f\|_{L^2}+\||x|^2f\|_{L^2}+\||v|^2f\|_{L^2}+\|Pf\|_{L^2}+\||x|^2Pf\|_{L^2}+\||v|^2Pf\|_{L^2}.$

To shorten the proof we shall consider here only the case
\[
  P = \theta + \nu v\cdot\nabla_x + \mu x\cdot\nabla_v + \beta v\cdot\nabla_v + \alpha\Delta_x + \sigma\Delta_v + \gamma\diver_v\nabla_x
\]
(cf. the definition of the operator $A$ in \fer{def:A}), but exactly the same
strategy extends to the case, where $P$ is a general quadratic polynomial.

First we define the mollifying delta sequence
\[
  \phi_n(x,v):= n^6\phi(nx,nv),\qquad n\in\N, \;x,v\in\R^3,
\]
where
\begin{eqnarray*}
  && \phi\in C_0^\infty(\R^6), \quad \phi(x,v)\ge 0, \\
  && \int\!\!\!\!\int\phi(x,v)dxdv = 1, \quad\mbox{and}\quad \supp\,\phi \subset \{|x|^2+|v|^2\le 1\}.
\end{eqnarray*}
By definition we have the following properties:
\begin{itemize}
\item[(I)] $\phi_n\to\delta$ in ${\cal D}'(\R^6),$
\item[(II)] $\frac{1}{n}\partial_{x_i}\phi_n,\,\frac{1}{n}\partial_{v_i}\phi_n \rightarrow 0$ in ${\cal D}'(\R^6),$ $i=1,2,3,$
\item[(III)] $(x,v)^\alpha \partial^{\beta}_{(x,v)} \left[ (x,v)^\gamma \phi_n(x,v) \right] \rightarrow 0 $ in ${\cal D}'(\R^6),$ with $\alpha,\beta,\gamma\in \N^{6}_{0}$ multi-indexes \\[0.17cm]
and $|\gamma|>0,$ since $(x,v)^\gamma\phi_n \rightarrow 0$ in ${\cal D}'(\R^6).$
\end{itemize}

The cutoff sequence is
\[
  \psi_n(x,v):=\psi\left( \frac{|(x,v)|}{n}\right),\qquad n\in\N, \;x,v\in\R^3,
\]
where $\psi$ satisfies
$$
  \psi\in C_0^\infty(\R), \quad  0\le\psi(z)\le 1, \quad \supp\psi \subset [-1,1], \quad \psi|_{[-\frac{1}{2},\frac{1}{2}]} \equiv 1,
$$
and
$$
  |\psi^{(j)}(z)|\le C_j, \quad \forall z\in\R, \;j=1,2 .
$$
The sequence $\psi_n$ has the following properties:
\begin{itemize}
\item[(IV)] $\psi_n \to 1$ pointwise,
\item[(V)]
$(x,v)^\alpha\partial_{(x,v)}^{\beta}\psi_n(x,v) = \frac{1}{n}\frac{(x,v)^\alpha (x,v)^\beta}{|(x,v)|}\psi'\left(\frac{|(x,v)|}{n}\right),$ with $\alpha, \beta \in \N^{6}_{0},$\\[0.1cm]
 $|\alpha| = |\beta| = 1,$ are supported in the annulus
\begin{center}
$\supp\left(\psi'\left(\frac{|(x,v)|}{n}\right)\right) \;=\; \left\{(x,v) \;|\; n/2 \leq |(x,v)| \leq n \right\} \;=:\; V_n,$
\end{center}
and they are in $L^{\infty}(\R^6),$ uniformly in $n\in\N.$
\item[(VI)]
$ n\partial_{(x,v)}^{\alpha}\psi_n(x,v) = \frac{(x,v)^\alpha }{|(x,v)|}\psi'\left(\frac{|(x,v)|}{n}\right),$ with $\alpha \in \N^{6}_{0},$ $|\alpha| = 1,$ \\[0.2cm]
are uniformly bounded in $L^{\infty}(\R^6).$
\item[(VII)] $\partial_{(x,v)}^{\alpha}\psi_n(x,v) = \frac{(x,v)^{\alpha}}{n^2|(x,v)|^2}\psi''\left(\frac{|(x,v)|}{n}\right) + \left(\frac{1}{n^2|(x,v)|} - \frac{(x,v)^{\alpha}}{n^3|(x,v)|^3} \right)\psi'\left(\frac{|(x,v)|}{n}\right),$ \\[0.2cm]
with  $|\alpha| = 2$ have support on $V_n$ and converge uniformly to $0$ in $L^{\infty}(\R^6).$
\end{itemize}

\medskip
We now define the approximating sequence
\[
  f_n(x,v) := (f\ast\phi_n)(x,v)\cdot\psi_n(x,v),\quad n\in\N,
\]
where `$\ast$' denotes the convolution in $x$ and $v.$ \\
By construction we have $f_n\in C^{\infty}_{0}(\R^6) = D(P).$

Since we can split our operator as
\begin{eqnarray*}
 P & = & \sum_{i=1}^3 \left[ \frac{\theta}{3} + \nu v_i\partial_{x_i} + \mu x_i\partial_{v_i} + \beta v_i\partial_{v_i} + \alpha\partial^{2}_{x_i} + \sigma\partial^{2}_{v_i} + \gamma\partial_{v_i}\partial_{x_i}\right] \\
& = & \sum_{i=1}^3 \tilde{p}(x_i,v_i,\partial_{x_i},\partial_{v_i}),
\end{eqnarray*}
we shall in the sequel only consider
$$
\tilde{P} \;=\; \tilde{p}(y,z,\partial_y,\partial_z), \quad y,z\in\R
$$
acting in one spatial direction $y=x_j$ and one velocity direction $z=v_j.$ \\
We have to prove that $f_n(x,v)\rightarrow f(x,v)$ in the graph norm
$$
\|f\|_{\tilde{P}}=\|f\|_{L^2}+\||x|^2f\|_{L^2}+\||v|^2f\|_{L^2}+\|\tilde{P}f\|_{L^2}+\||x|^2\tilde{P}f\|_{L^2}+\||v|^2\tilde{P}f\|_{L^2}.
$$

According to the 6 terms of the graph norm we split the proof into 6 steps:

{\bf Step 1: }
By applying (P1) and (P4), we have
\[
  f_n\to f \quad \mbox{in} \quad L^2(\R^6).
\]\\[-0.5cm]

{\bf Step 2: }
For the second term of the graph norm we write
\[
  x_{i}^{2} f_n = (x_{i}^{2}f\ast\phi_n)\psi_n + 2( x_{i}f\ast  x_{i}\phi_n)\psi_n + (f\ast x_{i}^{2}\phi_n)\psi_n.
\]
The first summand converges to $ x_{i}^{2}f$ in $L^2(\R^6)$ and both the second and the third terms
converge to $0$ by (III), since also $x_if$ belongs to $L^2(\R^6)$ by interpolation.\\[0.2cm]

{\bf Step 3: }
For the third term of the graph norm the same argument as in previous step can be used.
Hence we have
\[
  f_n\to f\qquad\mbox{in}\quad Z.
\]\\[-0.5cm]

{\bf Step 4: }
To prove that $\tilde{P}f_n\to \tilde{P}f$ in $L^2(\R^6)$ we write:
\begin{eqnarray*}
  \tilde{P}f_n & = & \frac{\theta}{3}(f\ast\phi_n)\psi_n +
  \nu(zf_y\ast\phi_n)\psi_n + \mu(yf_z\ast\phi_n)\psi_n + \beta(zf_z\ast\phi_n)\psi_n\\
       &   & +\alpha(f_{yy}\ast\phi_n)\psi_n + \sigma(f_{zz}\ast\phi_n)\psi_n +
             \gamma(f_{yz}\ast\phi_n)\psi_n +r_n^1(y,z)\\
       & = & (\tilde{P} f\ast\phi_n)\psi_n + r_n^1(y,z).
\end{eqnarray*}
As we shall show, all thirteen terms of the remainder
\begin{eqnarray*}
  r_n^1\! & = & \nu(f\ast \partial_y(z\phi_n))\psi_n + \nu(f\ast \phi_n)z\partial_y\psi_n +
  \mu(f\ast y\partial_z\phi_n)\psi_n \\
& & + \mu(f\ast \phi_n)y\partial_z\psi_n + \beta(f\ast\partial_z(z\phi_n))\psi_n + \beta(f\ast\phi_n))z\partial_z\psi_n \\
& & + 2\alpha(f\ast(\frac{1}{n}\partial_y\phi_n))(n\partial_y\psi_n) + \alpha(f\ast\phi_n))(\partial_{y}^{2}\psi_n) + 2\sigma(f\ast\frac{1}{n}\partial_z\phi_n)n\partial_z\psi_n \\
& & + \sigma(f\ast\phi_n)\partial_{z}^{2}\psi_n + \gamma(f\ast(\frac{1}{n}\partial_z\phi_n))(n\partial_y\psi_n) +  \gamma(f\ast(\frac{1}{n}\partial_y\phi_n))(n\partial_z\psi_n) \\
& & + \gamma(f\ast \phi_n)\partial_{y}\partial_z\psi_n
\end{eqnarray*}
converge to $0$ in $L^2(\R^6).$\\
The first, the third and the fifth terms converge to $0$ in $L^2(\R^6)$ by (III).\\
In the second, fourth and the sixth terms, exploiting (V) we have

\begin{equation}
\label{VIII}
 \Vert (f\ast\phi_n)(z\partial_y\psi_n)\Vert_{L^2(R^6)} \leq C\Vert f\ast\phi_n - f\Vert_{L^2(V_n)} +
\Vert f \Vert_{L^2(V_n)} \rightarrow 0,
\end{equation}

because $\Vert f\Vert_{L^2(R^6)} = \Vert f\Vert_{L^2(B_{1/2}(0))} + \sum_{k=0}^{\infty} \Vert f\Vert_{L^2(V_{2^k})}.$\\
For what the seventh, ninth, eleventh and twelfth terms are concerned, we can exploit (VI) and then (II).\\
The remaining terms can be handled thanks to (VII). \\[0.2cm]

{\bf Step 5: }
To prove that $|x|^2\tilde{P}f_n\to |x|^2Pf$ in $L^2(\R^6)$ we write:
\begin{eqnarray*}
  x_{i}^{2}\tilde{P}f_n & = & \frac{\theta}{3}(x_{i}^{2}f\ast\phi_n)\psi_n +
  \nu(x_{i}^{2}zf_y\ast\phi_n)\psi_n + \mu(x_{i}^{2}yf_z\ast\phi_n)\psi_n +
  \beta(x_{i}^{2}zf_z\ast\phi_n)\psi_n \\
& & +\alpha(x_{i}^{2}f_{yy}\ast\phi_n)\psi_n + \sigma(x_{i}^{2}f_{zz}\ast\phi_n)\psi_n +
   \gamma(x_{i}^{2}f_{yz}\ast\phi_n)\psi_n +r_n^2(y,z)\\
& = & (x_{i}^{2}\tilde{P}f\ast\phi_n)\psi_n + r_n^2(y,z).
\end{eqnarray*}

The remainder $r_n^2$ can be split in the following way ($y=x_j,$ $z=v_j$):
\begin{eqnarray*}
r^{2}_{n,\theta}\! & = &  \frac{2}{3}\theta(x_{i}f\ast y\phi_n)\psi_n + \frac{\theta}{3}(f\ast x_{i}^{2}\phi_n)\psi_n \\
r^{2}_{n,\nu}\! & = & 2\nu(zx_{i}f\ast \partial_y(x_{i}\phi_n))\psi_n -  2\nu\delta_{ij}(zf\ast x_{i}\phi_n)\psi_n + \nu(zf\ast \partial_y(x_{i}^{2}\phi_n))\psi_n \\
& & + \nu(x_{i}^{2}f\ast z\partial_y\phi_n)\psi_n  + 2\nu(x_{i}f\ast x_{i}z\partial_y\phi_n)\psi_n + \nu(f\ast x_{i}^{2}z\partial_y\phi_n)\psi_n \\
& & + \nu(x_{i}^{2}f\ast \phi_n)z\partial_y\psi_n + 2\nu(x_{i}f\ast x_{i}\phi_n)z\partial_y\psi_n +  \nu(f\ast x_{i}^{2}\phi_n)z\partial_y\psi_n \\
r^{2}_{n,\mu}\! & = &  2\mu(x_{i}yf\ast \partial_z(x_{i}\phi_n))\psi_n + \mu(yf\ast \partial_z(x_{i}^{2}\phi_n))\psi_n + \mu(x_{i}^{2}f\ast y\partial_z\phi_n)\psi_n \\
& & + 2\mu(x_{i}f\ast x_{i}y\partial_z\phi_n)\psi_n + \mu(f\ast x_{i}^{2}y\partial_z \phi_n)\psi_n + \mu(x_{i}^{2}f\ast \phi_n)y\partial_z\psi_n \\
& & + 2\mu(x_{i}f\ast x_{i}\phi_n)y\partial_z\psi_n + \mu(f\ast x_{i}^{2}\phi_n)y\partial_z\psi_n \\
r^{2}_{n,\beta}\! & = &  2\beta(x_{i}zf \ast x_{i}\partial_z\phi_n)\psi_n - 2\beta(x_{i}f\ast x_{i}\phi_n)\psi_n +
 \beta(zf\ast x_{i}^{2}\partial_z\phi_n)\psi_n - \beta(f\ast x_{i}^{2}\phi_n)\psi_n \\
& & +\beta(x_{i}^{2}f\ast \partial_z(z\phi_n))\psi_n +  2\beta(x_{i}f\ast x_{i}\partial_z(z\phi_n))\psi_n + \beta(f\ast x_{i}^{2}\partial_z(z\phi_n))\psi_n \\
& & +\beta(x_{i}^{2}f\ast \phi_n)z\partial_z\psi_n + 2\beta(x_{i}f\ast x_{i}\phi_n)z\partial_z\psi_n + \beta(f\ast x_{i}^{2}\phi_n)z\partial_z\psi_n   \\
r^{2}_{n,\alpha}\! & = & 2\alpha(x_{i}f \ast \partial_{y}^{2}(x_{i}\phi_n))\psi_n - 4\alpha\delta_{ij}(f \ast \partial_{y}(x_{i}\phi_n))\psi_n + \alpha(f \ast \partial_{y}^{2}(x_{i}^{2}\phi_n))\psi_n \\
& & +2\alpha(x_{i}^{2}f \ast \frac{1}{n}\partial_{y}\phi_n)n\partial_y\psi_n + 4\alpha(x_{i}f \ast \frac{x_{i}}{n}\partial_{y}\phi_n)n\partial_y\psi_n  + 2\alpha(f \ast \frac{x_{i}^{2}}{n}\partial_{y}\phi_n)n\partial_y\psi_n \\
& & +\alpha(x_{i}^{2}f \ast \phi_n)\partial_{y}^{2}\psi_n + 2\alpha(x_{i}f \ast x_{i}\phi_n)\partial_{y}^{2}\psi_n + \alpha(f \ast x_{i}^{2}\phi_n)\partial_{y}^{2}\psi_n \\
r^{2}_{n,\sigma}\! & = & 2\sigma(x_{i}f \ast x_{i}\partial_{z}^{2}\phi_n)\psi_n + \sigma(f \ast x_{i}^{2}\partial_{z}^{2}\phi_n)\psi_n + 2\sigma(x_{i}^{2}f \ast \frac{1}{n}\partial_{z}\phi_n)n\partial_z\psi_n \\
& & + 4\sigma(x_{i}f \ast \frac{x_{i}}{n}\partial_{z}\phi_n)n\partial_z\psi_n + 2\sigma(f \ast \frac{x_{i}^{2}}{n}\partial_{z}\phi_n)n\partial_z\psi_n + \sigma(x_{i}^{2}f \ast \phi_n)\partial_{z}^{2}\psi_n \\
& & + 2\sigma(x_{i}f \ast x_{i}\phi_n)\partial_{z}^{2}\psi_n + \sigma(f \ast x_{i}^{2}\phi_n)\partial_{z}^{2}\psi_n  \\
r^{2}_{n,\gamma}\! & = & 2\gamma (x_{i}f \ast \partial_{y}(x_{i}\partial_z\phi_n))\psi_n - 2\gamma\delta_{ij}(f \ast x_{i}\partial_{z}\phi_n)\psi_n + \gamma (f \ast \partial_{y}\partial_z(x_{i}^{2}\phi_n))\psi_n \\
& & + \gamma(x_{i}^{2}f \ast \frac{1}{n}\partial_{z}\phi_n)n\partial_y\psi_n + 2\gamma(x_{i}f \ast \frac{x_{i}}{n}\partial_{z}\phi_n)n\partial_y\psi_n + \gamma(f \ast \frac{x_{i}^{2}}{n}\partial_{z}\phi_n)n\partial_y\psi_n \\
& & + \gamma(x_{i}^{2}f \ast \frac{1}{n}\partial_{y}\phi_n)n\partial_z\psi_n + 2\gamma(x_{i}f \ast \frac{x_{i}}{n}\partial_{y}\phi_n)n\partial_z\psi_n + \gamma(f \ast \frac{x_{i}^{2}}{n}\partial_{y}\phi_n)n\partial_z\psi_n \\
& & +\gamma(x_{i}^{2}f \ast \phi_n)\partial_y\partial_z\psi_n + 2\gamma(x_{i}f \ast x_{i}\phi_n)\partial_y\partial_z\psi_n + \gamma(f \ast x_{i}^{2}\phi_n)\partial_y\partial_z\psi_n.
\end{eqnarray*}
By the properties (I)-(VII) and estimate like (\ref{VIII}), it can be easily seen that each term converges to $0$ in $L^2(\R^6).$\\[0.2cm]

{\bf Step 6: }
In analogy to $|x|^2\tilde{P}f_n,$ the sequence $|v|^2\tilde{P}f$ can be split as
$$
  v_{i}^{2}\tilde{P}f_n = (v_{i}^{2}\tilde{P}f\ast\phi_n )\psi_n + r_{n}^{3}(y,z).
$$
Due to the symmetry of the operator $\tilde{P}$ in $x$ and $v,$ the terms of the remainder $r_{n}^{3}$ can be obtained from $r_{n}^{2}$ by interchanging $y$ and $z$ (and changing the coefficients), except for the following term
\begin{eqnarray*}
  v_{i}^{2}[\beta z\partial_z\left((f\ast\phi_n)\psi_n\right)] & = & \beta(v_{i}^{2}zf_z\ast\phi_n)\psi_n + r^{3}_{n,\beta},
\end{eqnarray*}
where
\begin{eqnarray*}
  r^{3}_{n,\beta}\! & = & 2\beta(v_{i}zf \ast \partial_z(v_{i}\phi_n))\psi_n - 2\beta(1+\delta_{ij})(v_{i}f \ast v_{i}\phi_n)\psi_n + \beta(zf \ast \partial_z(v_{i}^{2}\phi_n))\psi_n \\
& & - \beta(f \ast v_{i}^{2}\phi_n)\psi_n + \beta(v_{i}^{2}f \ast \partial_z(z\phi_n))\psi_n + 2\beta(v_{i}f \ast v_{i}\partial_z(z\phi_n))\psi_n \\
&& + \beta(f \ast v_{i}^{2}\partial_z(z\phi_n))\psi_n + \beta(v_{i}^{2}f \ast \phi_n)z\partial_z\psi_n + 2\beta(v_{i}f \ast v_{i}\phi_n)z\partial_z\psi_n \\
& &+ \beta(f \ast v_{i}^{2}\phi_n)z\partial_z\psi_n
\end{eqnarray*}
converges to $0$ in $L^2(\R^6),$ since (I)-(VII) and (\ref{VIII}) can be used.\qed

\subsection*{Proof of Proposition \ref{prop:decay}}

First, we shall prove the following estimates on the
derivatives of the Green's function (\ref{eq:Green}):
\begin{eqnarray}
\label{eq:v-decayOfGreen} |\nabla_v G(t,x,v,x_0,v_0)| & \leq & b
\frac{G(t,\frac{x}{2},\frac{v}{2},\frac{x_0}{2},\frac{v_0}{2})}{\sqrt
t},\quad \forall \,t\leq t_0,\\ \label{eq:x-decayOfGreen} |\nabla_x
G(t,x,v,x_0,v_0)| & \leq &
b'\frac{G(t,\frac{x}{2},\frac{v}{2},\frac{x_0}{2},\frac{v_0}{2})}{\sqrt
t},\quad \forall \,t \leq t_1.
\end{eqnarray}
with $b=b(\alpha,\gamma, \sigma), t_0=t_0(\alpha,\beta,\sigma,\gamma), b'=b'(\alpha,\gamma,\sigma)$ and
$t_1=t_1(\alpha,\beta,\sigma,\gamma).$ The $v$-derivative of $G$ is given by
\begin{eqnarray}
\label{eq:v-derOfGreen}
\nabla_v G(t,x,v,x_0,v_0)&=&G(t,x,v,x_0,v_0)\left[- \frac{\left(\mu(t)e^{\beta t} - 2\nu(t)\frac{e^{\beta t}-1}{\beta}\right)(x-\frac{e^{\beta t}-1}{\beta}v-x_0)}{f(t)} \right. \nonumber\\
&& -\left. \frac{\left(2\lambda(t)e^{\beta t} - \mu(t)\frac{e^{\beta t}-1}{\beta}\right)(e^{\beta t}v-v_0)}{f(t)}\right].
\end{eqnarray}
For all real $a,b,c > 0$ such that $c/\sqrt{a}\leq b\sqrt{2e},$ one easily verifies that
\begin{equation}
\label{eq:ineq}
  c|x| \;\:\leq\;\:b e^{a|x|^2}, \quad\forall\, x\in\R^3.
\end{equation}

Since $\alpha,\sigma>0,$ we have for $t>0$ small enough
$$
\nu(t)-\frac{1}{2}\mu(t)\;\:>\;\: 0, \qquad  \lambda(t)-\frac{1}{2}\mu(t)\;\:>\;\: 0.
$$
In order to apply the estimate (\ref{eq:ineq}) to the two terms inside the squared bracket in (\ref{eq:v-derOfGreen})
we shall use for $t$ small:
$$
\frac{c_1}{\sqrt{a_1}}\;\::=\:\;\frac{\frac{\sqrt{t}}{f(t)}\left|\mu(t)e^{\beta t}-2\nu(t)\frac{e^{\beta
t}-1}{\beta}\right|}{\sqrt{\frac{3}{4}\frac{\nu(t)-\frac{1}{2}\mu(t)}{f(t)}}}\;\:\sim\:\; \frac{2\gamma}{\sqrt{3(\alpha \sigma-\gamma^2)(\sigma+\gamma)}}\;\:\leq\:\; b_1\sqrt{2e},$$
with $b_1=\gamma/\sqrt{3(\alpha \sigma-\gamma^2)(\sigma+\gamma)}.$
Similarly,
$$
\frac{c_2}{\sqrt{a_2}}\;\::=\;\:\frac{\frac{\sqrt{t}}{f(t)}\left|2\lambda(t)e^{\beta t} -\mu(t)\frac{e^{\beta
t}-1}{\beta}\right|}{\sqrt{\frac{3}{4}\frac{\lambda(t)-\frac{1}{2}\mu(t)}{f(t)}}}\;\:\sim\;\: \frac{2\alpha}{\sqrt{3(\alpha \sigma-\gamma^2)(\alpha+\gamma)}}\;\:\leq\;\: b_2\sqrt{2e}, $$
with $b_2=\alpha/\sqrt{3(\alpha \sigma-\gamma^2)(\alpha+\gamma)}.$ Then, there exists some $t_0>0$ such that, for all
$t\leq t_0,$ the two inequalities can be combined with $b=\max\{b_1,b_2\}$ to give
$$
\!\!\!\!\!\!\!\!\!\!\!\!\!\!\!\!\!\left|\frac{\left(\mu(t)e^{\beta t} -
2\nu(t)\frac{e^{\beta
t}-1}{\beta}\right)\left(x-\frac{e^{\beta t}-1}{\beta}v-x_0\right) +
\left(2\lambda(t)e^{\beta t} - \mu(t)\frac{e^{\beta
t}-1}{\beta}\right)(e^{\beta t}v-v_0)}{f(t)}\right|\sqrt{t} \:\leq
$$
\vspace*{-0.5cm}
\begin{eqnarray*}
&&\!\!\!\!\!\!\!\!\!\!\!\leq \frac{\sqrt{t}}{f(t)}\left\{\left|\mu(t)e^{\beta t} - 2\nu(t)\frac{e^{\beta t}-1}{\beta}\right|\left|x-\frac{e^{\beta t}-1}{\beta}v-x_0\right| + \left|2\lambda(t)e^{\beta t} - \mu(t)\frac{e^{\beta t}-1}{\beta}\right||e^{\beta t}v-v_0|\right\} \\
&&\!\!\!\!\!\!\!\!\!\!\!\leq b \exp\left\{\frac{\left(\nu(t)-\frac{1}{2}\mu(t)\right)\left|x-\frac{e^{\beta t}-1}{\beta}v-x_0\right|^2 + \left(\lambda(t)-\frac{1}{2}\mu(t)\right)\left|e^{\beta t}v-v_0\right|^2 }{\frac{4}{3}f(t)}\right\}\\
&&\!\!\!\!\!\!\!\!\!\!\!\leq b\exp\left\{\frac{\nu(t)\left|x-\frac{e^{\beta t}-1}{\beta}v-x_0\right|^2 + \lambda(t)\left|e^{\beta t}v-v_0\right|^2 + \mu(t)\left(x-\frac{e^{\beta t}-1}{\beta}v-x_0\right)\cdot\left(e^{\beta t}v-v_0\right)}{\frac{4}{3}f(t)}\right\}\!.
\end{eqnarray*}
Hence,
\begin{eqnarray*}
&&\!\!\!\!\!\!\!\!\!\!\!|\nabla_v G(t,x,v,x_0,v_0)| \;\:\leq\;\: b \frac{G(t,x,v,x_0,v_0)}{\sqrt{t}}\\
&&\!\!\!\!\!\!\!\!\!\!\!\times\exp\left\{\frac{3}{4}\frac{\nu(t)\left|x-\frac{e^{\beta t}-1}{\beta}v-x_0\right|^2 + \lambda(t)\left|e^{\beta t}v-v_0\right|^2 + \mu(t)\left(x-\frac{e^{\beta t}-1}{\beta}v-x_0\right)\cdot\left(e^{\beta t}v-v_0\right)}{f(t)}\right\}
\end{eqnarray*}
and the decay (\ref{eq:v-decayOfGreen}) follows by comparison with (\ref{eq:Green}).\\

Next we consider the $x$-derivative of the Green's function,
\begin{eqnarray*}
\nabla_x G(t,x,v,x_0,v_0) & = & G(t,x,v,x_0,v_0)\left[- \frac{2\nu(t)(x-\frac{(e^{\beta t}-1)}{\beta}v-x_0) + \mu(t)(e^{\beta t}v-v_0)}{f(t)}\right].
\end{eqnarray*}
Analogously, the decay (\ref{eq:x-decayOfGreen}) follows by exploiting that for $t$ small enough\\
$$ \frac{\frac{\sqrt{t}}{f(t)}
2\nu(t)}{\sqrt{\frac{3}{4}\frac{\nu(t)-\frac{1}{2}\mu(t)}{f(t)}}} \;\:\sim\;\:\frac{2\sigma}{\sqrt{3(\alpha \sigma-\gamma^2)(\sigma+\gamma)}}\;\:\leq\;\: b_1'\sqrt{2e},$$
$$\frac{\frac{\sqrt{t}}{f(t)}|\mu(t)|}{\sqrt{\frac{3}{4}\frac{\lambda(t)-\frac{1}{2}\mu(t)}{f(t)}}}\;\:\sim\;\:\frac{2\gamma}{\sqrt{3(\alpha \sigma-\gamma^2)(\alpha+\gamma)}}\;\:\leq\;\: b_2'\sqrt{2e},$$
with appropriate $b_1'(\alpha,\gamma, \sigma), b_2'(\alpha,\gamma, \sigma).$ \\

Since
\begin{eqnarray*}
e^{t\overline{A}}w_0(x,v) & = & \int\!\!\!\int\!
G(t,x,v,x_0,v_0)w_0(x_0,v_0)\,dx_0\,d v_0,
\end{eqnarray*}
we have
\begin{eqnarray}
\label{eq:decay}
|\nabla_v e^{t\overline{A}}w_0(x,v)| & \leq & \int\!\!\!\int\!
|\nabla_v G(t,x,v,x_0,v_0)|\,|w_0(x_0,v_0)| \,dx_0\,d v_0 \nonumber \\
& \leq &
b t^{-1/2} \int\!\!\!\int\!
G\left(t,\frac{x}{2},\frac{v}{2},\frac{x_0}{2},\frac{v_0}{2}\right)|w_0(x_0,v_0)|
\,dx_0\,d v_0 \nonumber\\
& = & 64 b\, t^{-1/2} \int\!\!\!\int\!
G(t,\tilde{x},\tilde{v},\tilde{x_0},\tilde{v_0})|w_0(2\tilde{x_0},2\tilde{v_0})|\,d\tilde{x_0}\,d\tilde{v_0} \nonumber \\
& = & 64 b\,t^{-1/2}
e^{t\overline{A}}\tilde{w} _0(\tilde{x},\tilde{v}),\quad \forall\,t\leq t_0.
\end{eqnarray}
Here we used the decay (\ref{eq:v-decayOfGreen}), and we put
$\tilde{x} = \frac{x}{2},$ $\tilde{v} = \frac{v}{2}$ and
$\tilde{w} _0(\tilde{x},\tilde{v}) = |w_0(2\tilde{x},2\tilde{v})|.$
The assertion (\ref{eq:decayII-v}) follows directly by applying
the estimate (\ref{eq:decayI}) to (\ref{eq:decay})
and choosing $T_0=\min\{t_0,t_1\}.$
\\ The estimate (\ref{eq:decayII-x}) can be obtained analogously.
\qed



\newpage
\textbf{Anton Arnold}\\
Institut f\"ur Numerische und Angewandte Mathematik \\
University of M\"unster \\
Einsteinstr. 62\\
D-48149 M\"unster, Germany \\

Email: aarnold@math.uni-muenster.de\\[30mm]

\textbf{Elidon Dhamo}\\
Institut f\"ur Numerische und Angewandte Mathematik \\
University of M\"unster \\
Einsteinstr. 62\\
D-48149 M\"unster, Germany \\

Email: dhamo@math.uni-muenster.de\\[30mm]

\textbf{Chiara Manzini}\\
Scuola Normale Superiore \\
Piazza dei Cavalieri 7\\
I-56126 Pisa, Italy\\

Email: cmanzini@sns.it

\end{document}